%% file: main_bit.tex
\RequirePackage{fix-cm}

\documentclass[smallextended,numbook,runningheads]{svjour3}     
\smartqed  
\usepackage{graphicx}
\usepackage{amsmath}
\usepackage{epstopdf} 
\usepackage{mathptmx}      
\usepackage[T1]{fontenc}
\usepackage[utf8]{inputenc}
\usepackage{newtxtext,newtxmath}
\usepackage{hyperref}
\hypersetup{breaklinks,colorlinks,linkcolor=blue,citecolor=blue}
\usepackage{cleveref}
\usepackage{amsmath}
\usepackage{tikz}
\usepackage{pgfplots}
\pgfplotsset{compat=1.17}
\usepackage{comment}
\usepackage{pierson}

\usetikzlibrary{shapes.geometric}

\journalname{BIT}

\begin{document}

\title{Parallel Scaling of the Regionally-Implicit Discontinuous Galerkin Method with Quasi-Quadrature-Free Matrix Assembly
\thanks{ This work was performed under the auspices of the U.S. Department of Energy by Lawrence Livermore National Laboratory under Contract DE-AC52-07NA27344.
Scaling studies were performed at Institute for Cyber-Enabled Research at Michigan State University and at Livermore Computing at Lawrence Livermore National Laboratory.
This work was funded in part by ONR High Order Scalable Solvers grant N0014-19-1-2476, and AFOSR Computational Non-Ideal Plasma Physics grant FA9550-17-1-0394, and DoE SciDAC TEAMS grant DE-SC0017955.
JAR was supported in part by NSF Grants DMS--1620128 and DMS--2012699.
}}



\author{
         Andrew J. Christlieb \and 
         Pierson T. Guthrey\footnote{Corresponding Author} \and 
         James A. Rossmanith  
}


\titlerunning{Parallel Scaling of the Regionally-Implicit Discontinuous Galerkin Method}
\authorrunning{A.J. Christlieb, P.T. Guthrey, and J.A. Rossmanith}

\institute{
 Andrew J. Christlieb \at 
 Michigan State University \\
 Department of Computational Mathematics, Science and Engineering \\
  428 S. Shaw Lane \\
  East Lansing, Michigan 48824, USA\\
 \email{\href{mailto:christli@msu.edu}{christli@msu.edu}}
 \and
 Pierson T. Guthrey \at
              Weapons and Complex Integration \\
              Lawrence Livermore National Laboratory\\ 
              Livermore, CA 94550, USA \\
              Tel.: +337-781-5574 \\
              \email{\href{mailto:guthrey1@llnl.gov}{guthrey1@llnl.gov}}
           \and
  James A. Rossmanith \at
  Iowa State University \\
 Department of  Mathematics \\
 411 Morrill Road \\
 Ames, Iowa 50011, USA \\
 \email{\href{mailto:rossmani@iastate.edu}{rossmani@iastate.edu}}
 }

\date{Received: date / Accepted: date}

\maketitle

\begin{abstract}
\input{secabstract.tex}
\end{abstract}


\input{secintroduction.tex}
\input{secridg.tex}
\input{secJacobians.tex}
\input{secterminology.tex}
\input{sec_implement.tex}
\input{section_1d.tex}

\input{section_2d.tex}

\input{section_3D_adv.tex}
\input{secfinal.tex}


\begin{acknowledgements}
This work was performed under the auspices of the U.S. Department of Energy by Lawrence Livermore National Laboratory under Contract DE-AC52-07NA27344.
The authors would like to thank Philipp Grete and Forrest Glines for their helpful conversations.
\end{acknowledgements}

\section*{Disclaimer} 
This document was prepared as an account of work sponsored by an agency of the United States government. Neither the United States government nor Lawrence Livermore National Security, LLC, nor any of their employees makes any warranty, expressed or implied, or assumes any legal liability or responsibility for the accuracy, completeness, or usefulness of any information, apparatus, product, or process disclosed, or represents that its use would not infringe privately owned rights. Reference herein to any specific commercial product, process, or service by trade name, trademark, manufacturer, or otherwise does not necessarily constitute or imply its endorsement, recommendation, or favoring by the United States government or Lawrence Livermore National Security, LLC. The views and opinions of authors expressed herein do not necessarily state or reflect those of the United States government or Lawrence Livermore National Security, LLC, and shall not be used for advertising or product endorsement purposes.

\bibliographystyle{spmpsci}      

%
%

\end{document}

%% file: secabstract.tex

In this work we investigate the parallel scalability of the numerical method developed in Guthrey and Rossmanith
[{\it The regionally implicit discontinuous Galerkin method: Improving the
stability of DG-FEM}, SIAM J. Numer. Anal. (2019)].
We develop an implementation of the regionally-implicit discontinuous Galerkin (RIDG) method in DoGPack, which is an open source C++ software package for discontinuous Galerkin methods. Specifically, we develop and test a hybrid OpenMP and MPI parallelized implementation of DoGPack with the goal of exploring the efficiency and scalability of RIDG in comparison to the popular strong stability-preserving Runge-Kutta discontinuous Galerkin (SSP-RKDG) method.  We demonstrate that RIDG methods are able to hide communication latency associated with distributed memory parallelism, due to the fact that almost all of the work involved in the method is highly localized to each element, producing a localized prediction for each region. We demonstrate the enhanced efficiency and scalability of the of the RIDG method and compare it to SSP-RKDG methods and show extensibility to very high order schemes.
The two-dimensional scaling study is performed on machines at the Institute for Cyber-Enabled Research at Michigan State University, using up to 1440 total cores on Intel(R) Xeon(R) Gold 6148 CPU @ 2.40GHz CPUs. The three-dimensional scaling study is
performed on Livermore Computing clusters at at Lawrence Livermore National Laboratory, using up to 28672 total cores on Intel Xeon CLX-8276L CPUs with Omni-Path interconnects.


%% file: secintroduction.tex

\section{Introduction}

In high-performance computing, {\it parallel scaling} is an important measure of algorithm efficiency, since it provides a quantifiable   
measure of how much computational benefit is gained by moving to more and more computing resources (i.e., more processors/cores/threads). There are two standard approaches measure scalability: (1) {\it weak scaling}, where the problem size (i.e., the degrees of freedom) increases with the number of cores such that the problem size per processor/core/thread is fixed; and (2) {\it strong scaling}, where the number of processors/cores/threads is increased while the problem size is fixed. For problems that are {\it CPU-bound}, i.e., problems where memory is not the primary concern, but computing the solution on a single processor takes a long time, the relevant scaling measure is strong scaling.

For the numerical solution of partial differential equations (PDEs), the standard strategy for strong scaling is to subdivide a problem of size $P$ over $n$ physical compute cores until a minimal problem size per physical compute core is reached (e.g., see Fischer, Heisey, and Min \cite{Fischer}). The ratio of problem size to number of compute cores is known as the granularity, $\eta$, and can be used to predict the ratio of MPI-communication and OpenMP overhead to parallelized work.  In the limit of minimal granularity $\eta \to 1$, parallelization overhead is maximized, reducing efficiency. 

In this work we are concerned with strong scaling via a hybrid OpenMP \cite{openmp} + MPI \cite{Gabriel2004} approach for a specific class of problems: hyperbolic conservation laws, with a specific type of spatial discretization: discontinuous Galerkin (DG) finite element methods (e.g., see Cockburn and Shu \cite{cockshu5}). The most common technique to time-advance DG spatial discretizations is the strong-stability-preserving Runge-Kutta (SSP-RK) scheme. 
This approach is explicit in time, and thus very efficient in terms of CPU and memory resources per time-step. However, a drawback is that SSP-RK schemes applied to high-order DG methods suffer from very small allowable time-steps, which then results in an overall reduction in efficiency (i.e., many time-steps to reach a fixed final time). In order to overcome these small time-step restrictions, Guthrey and Rossmanith \cite{RIDG_paper_2019} developed a new time-stepping approach for DG that was dubbed the {\it regionally-implicit DG} (RIDG) scheme. The primary goal of the current work is to develop, implement, and test a parallelized version of the regionally-implicit discontinuous Galerkin (RIDG) method.

\subsection{Scalability challenges}
The algorithms discussed in this paper are high-order discontinuous Galerkin (DG) methods with explicit time-stepping for hyperbolic conservation laws.  Explicit time-stepping algorithms for DG methods are usually associated with small time-steps, typically with Courant-Friedrichs-Lewy (CFL) numbers much less than $1$ (and indeed inversely proportional to the method order of accuracy).  The benefit of such schemes is that they have a very low computational cost associated with the (typically) nearest neighbor stenciling. Thus, there are two challenges associated with attempts to scale such schemes to multiple processors.  First, as we divide the computational domain into several pieces, each piece must synchronize with its neighbors across the domain-decomposition pseudo-boundaries, a process also known as {\it synchronization of the halo regions} \cite{Consortium2017}.  This must occur at least once per time-step (typically several times for multi-stage methods), and thus must be performed many times over the course of forming a solution for some fixed final time.  Secondly, for each time-step all processes must synchronize and compute the global maximum wave speed, so that each process is using the same time-step size given by the CFL-restriction.  This synchronization takes the form of an all-to-all reduction.  This can be potentially avoided by using local-time stepping, but we do not discuss this in this paper.

In the face of the challenges, discontinuous Galerkin methods enjoy enhanced scalability properties due to the nearest-neighbor stenciling typical of such methods and the highly element-localized compute intensity.
 Adaptive mesh refinement has been shown to work very well when applied to DG methods at scale \cite{Dumbser2013a}.
The arbitrary high order schemes using derivatives (ADER-DG) are predictor-corrector schemes that avoid repeated communication for each time-step by being fully discrete \cite{Dumbser2018}.  That is, there is a space-time basis underlying the time-step updates as opposed to a method-of-lines update of a purely spatial basis.  This approach shows improved MPI scaling efficiency compared to RKDG schemes, resulting in up to $23\%$ smaller $L_2$ errors with $1-2\%$ longer runtimes.  One goal of this paper is to perform a similar comparison of RIDG and RKDG methods in this paper. 

For discontinuous Galerkin methods, the granularity $\eta$ is usually discussed in terms of the total number of elements $P$, as opposed to the total number of problem degrees of freedom. This is because domain-decomposition is not typically implemented below the level of the elements.  That is, individual elements are not typically broken up as any part of a domain decomposition strategy.   

Domain decomposition for strong scaling such a problem involves dividing the physical domain into a number of sub-domains equal to the number of MPI tasks or compute nodes, each with a number of physical compute cores. This gives rise to domain decomposition pseudo-boundaries, in addition to the physical boundaries. In this paper we consider periodic boundary conditions, where physical boundaries become pseudo-boundaries, as they separate problem subdomains.  

    As discussed above, dynamic explicit time-stepping relies on an all-to-all reduction to ensure that each time-step taken uses a global estimate for the maximum wave speed in the computational domain.  That is, each subdomain assigned to an MPI task must compute its ``task-local'' maximum wave speed, and a ``task-global'' maximum wave speed must be computed via an all-to-all reduction to ensure that the Courant-Friedrichs-Lewy (CFL) condition is satisfied by all of the sub-domains. Since this procedure merely involves comparing wave speeds to find the maximum wave speed among all processes, the communication latency of this procedure is dominated by the total number of processes as $\log_2 P$, not the granularity $\eta$.  
    
    Every subdomain must perform halo-region communications with all subdomains with which it shares a pseudo-boundary.  This implies that periodic boundary conditions are more difficult to scale than Dirichlet or outflow conditions, since they inherently involve extra communications, specifically the {\it wrap-around} communications (e.g., the east side of the mesh communicates with the west side of the mesh).  As the granularity increases, the time spent copying data to MPI buffers {\it relative} to the work in the rest of the problem increases, potentially degrading efficiency.  However, this is a memory operation and thus can be relatively quick, even at low granularity. 
    
    Explicit time-stepping methods usually exhibit restrictive CFL conditions (even more-so with high-order DG spatial discretizations), which increase the total number of time steps, which in turn increases the total amount of communications per simulation time.   Thus, the scalability of an explicit time-stepping scheme depends on the scheme's ability to overcome communication overheads and hide latency as $\eta\to 1$. 
    
\subsection{Scope of this work}
    The purpose of the current work is to develop, implement, and study an OpenMP + MPI parallelization of the regionally implicit discontinuous Galerkin scheme (RIDG) \cite{Guthrey2017,RIDG_paper_2019}. In particular, we are concerned with comparing it to parallel implementations of the most commonly used variant of the discontinuous Galerkin method: the strong-stability-preserving Runge-Kutta DG (SSP-RKDG) scheme
    \cite{cockshu5,article:GoShu98,gottliebShuTadmor01}. 
    In  \cref{sec:ridg_review} we will briefly review the regionally-implicit method, and in particular, the two key parts of the time-stepping strategy: (1) regionally-implicit prediction and (2) explicit correction. 
In  \cref{sec:jacobians} we develop a quasi-quadrature-free approach for the Jacobian matrix assembly that is required in the regionally-implicit prediction step. Next we briefly define the terminology needed to quantify the accuracy and efficiency of our parallel implementations in  \cref{sec:terminology}.
A brief discussion and explanation of the proposed parallel implementation is provided in
\cref{sec:parallel}. 
In \cref{sec:1dresults}, \cref{sec:2dresults}, and \cref{sec:3dresults} we report our findings for examples in one, two, and three spatial dimensions, respectively.
The two-dimensional scaling study presented in \cref{sec:2dresults} is
performed on the Institute for Cyber-Enabled Research at Michigan State University's Intel(R) Xeon(R) Gold 6148 CPU @ 2.40GHz CPUs, 
 using up to 1440 total cores. The three-dimensional scaling study presented in 
 \cref{sec:3dresults} is
performed on Livermore Computing clusters at at Lawrence Livermore National Laboratory, using up to 28672 total cores on Intel Xeon CLX-8276L CPUs with Omni-Path interconnects. In these sections we demonstrate the benefits of using RIDG over SSP-RKDG, especially as the number of processors increase and as the order of the method increases. Finally, we conclude our findings in \cref{sec:conclusions}.
   

%% file: secridg.tex

\section{Regionally-implicit discontinuous Galerkin (RIDG) method}
\label{sec:ridg_review}
\input{tikz/ridg2d}
The regionally implicit discontinuous Galerkin (RIDG) method was developed in recent work by Guthrey and Rossmanith \cite{RIDG_paper_2019}. The goal of that work was to improve the linear stability of the the ADER-DG scheme \cite{article:Dumbser2006,article:GasDumHinMun2011,article:Zanotti2015}.
In particular, the origin of the RIDG method is tied to the method of Gassner et al. \cite{article:GasDumHinMun2011}, where it was shown that Lax-Wendroff discontinuous Galerkin \cite{Qiu2005a} can be formulated as a predictor-corrector method. The predictor is a local version of a spacetime DG method \cite{article:KlaVegVen2006,article:Sudirham2006}  (i.e., the predictor is something like a block-Jacobi update for a fully implicit spacetime DG method), and the corrector is an explicit method that uses the spacetime reconstructed solution from the predictor step. The drawback of the locally-implicit approach is that the resulting ADER-DG scheme suffers from small maximum allowable time-steps that decrease with increasing order of accuracy. The RIDG method instead uses a {\it regionally-implicit} prediction step, which has the benefit of resulting in a DG method with maximum allowable time-steps that are much closer to the optimal limit; and, in particular, the maximum allowed time-step does not degrade as the order of accuracy is increased. The drawback of RIDG is that per time-step it is more computational expensive than standard ADER-DG; however, this penalty is more than offset by the significantly increased stable time-step, especially for very high-order schemes. In this section we briefly review the key concepts of the RIDG approach. For more details we refer the reader to \cite{RIDG_paper_2019}.

\subsection{General setup}
Consider hyperbolic conservation laws of the form
\begin{equation}
\label{eqn:conslaw}
\vec{q}_{,t} + \div \mat{F}\left( \vec{q} \right) = \vec{0},
\end{equation}
where $\vec{q}\left(t,\vec{x} \right): \reals^+ \times \reals^\mdim \mapsto \reals^\meq$ is the vector of conserved variables, $\mat{F}\left(\vec{q}\right): \reals^\meq \mapsto \reals^{\meq \times \mdim}$ is the flux function, $\mdim$ is the number of spatial dimensions, and $\meq$ is the number of conserved variables.  We define $\Omega \subset \reals^\mdim$ to be a polygonal domain with boundary $\partial \Omega$, and discretize $\Omega$ using a finite set of non-overlapping elements, $\Tm_i$, such
that $\cup_{i=1}^\melems \Tm_i = \Omega$, where $\melems$ is the total number of elements.
Let ${\mathbb Q}\left(\mdeg, \mdim \right)$ denote the set of polynomials from $\reals^\mdim$ to $\reals$ 
with maximal polynomial degree $\mdeg$.
On the mesh of $\melems$ elements we define the {\it broken} finite element space:
\begin{equation*}
    \WS^h := \left\{ \vec{w}^h \in \left[ L^{\infty}(\Omega) \right]^{\meq}: \,
    \vec{w}^h \bigl|_{\Tm_i} \in \left[ {\mathbb Q} \left(\mdeg, \mdim \right) \right]^{\meq} \, \, \forall \Tm_i \right\},
\end{equation*}
where $h$ is the grid spacing.

Let $\Phi_{k}\left(\vec{x}\right)$ for $k=1,\ldots,\mbasis$ be a basis that spans 
${\mathbb Q} \left(\mdeg, \mdim \right)$ over $\Tm_i$ and is orthonormal:
\begin{equation*}
\frac{1}{|\Tm_i|} \int_{\Tm_i} \Phi_{k}\left(\vec{x}\right) 
\Phi_{\ell}\left(\vec{x}\right) \, d\vec{x} = \delta_{k\ell},
\end{equation*}
where $\delta_{k\ell}$ is the Kronecker delta and $|\Tm_i|$ is the volume
of element $\Tm_i$.

In order to setup the RIDG method we assume that on each element the
solution is of the following form:
\begin{equation}
\label{eqn:phi_basis_1d}
\vec{q}^h\left(t^n, \vec{x} \right) \Bigl|_{\Tm_i} = \vec{\Phi}\left(\vec{x}\right)^T \mat{Q}_i = \sum_{\ell=1}^{\mbasis} \vec{Q}_{i}^{\ell}\left(t^n\right) \, \Phi_{\ell}\left(\vec{x}\right),
\end{equation}
where $\vec{Q}_{i}^{\ell}(t^n)$ are the unknown degrees of freedom at
time level $t^n$.
Using this ansatz,  multiplying \cref{eqn:conslaw} by 
$\varphi_{k}$, integrating over the element $\Tm_i$, using integration-by-parts in space, and integrating over the
time slab $\left[t^n, t^{n+1} \right]$, results in the
following equation:
\begin{equation}
\label{eqn:semi_discrete_dg}
\begin{split}
\vec{Q}_{i}^{k}\left( t^{n+1} \right) &= 
\vec{Q}_{i}^{k}\left( t^{n} \right)
+ \frac{1}{|\Tm_i|} \int_{t^n}^{t^{n+1}} \int_{\Tm_i} 
\grad \Phi_{k} \cdot \mat{F}\left( \, \vec{q}^h \right)  \, d\vec{x} \, dt \\
&- \frac{1}{|\Tm_i|}\int_{t^n}^{t^{n+1}} \oint_{\partial \Tm_i} \Phi_k \, \vec{{\mathcal F}}\left(\vec{q}^h_{+}, \vec{q}^h_{-};
\vec{n} \right) \, d\vec{s} \, dt,
\end{split}
\end{equation}
where $\vec{n}$ is an outward-pointing normal vector to $\partial \Tm_i$, $\vec{q}^h_{+}$ and $\vec{q}^h_{-}$ are the states on either side of the boundary $\partial \Tm_i$, and $\vec{{\mathcal F}}$ is some appropriate numerical flux. Equation \cref{eqn:semi_discrete_dg} has the look and feel of a fully-discrete numerical method; however, the above volume and surface integral cannot be evaluated directly since they require knowledge of the solution over the entire time-slab: $\left[t^n, t^{n+1} \right]$. 

The idea of schemes such as ADER-DG \cite{article:Dumbser2006,article:GasDumHinMun2011,article:Zanotti2015},
its close cousin Lax-Wendroff DG \cite{Qiu2005a}, and indeed RIDG
\cite{RIDG_paper_2019}, is to form a {\it prediction} for solution 
$\vec{q}^h$ over the time-slab: $\left[t^n, t^{n+1} \right]$, and then to insert this prediction into equation \cref{eqn:semi_discrete_dg} to obtain a fully discrete update formula. The main difference between ADER-DG, Lax-Wendroff DG, and RIDG is only in how the prediction step is formed. We briefly review the details for RIDG below.

\subsection{RIDG prediction step}
\label{sec:RIDG1D_predict}
In order to describe the prediction step, we define {\it spacetime elements}
and {\it regions}:
\begin{equation}
\label{eqn:region}
{\mathcal S}_i := \left[ t^n, \, t^{n+1} \right] \times \Tm_{i} \qquad \text{and}
\qquad 
{\mathcal R}_i :=  \bigcup_{j: {\mathcal T}_i 
 \cap  {\mathcal T}_j \ne \emptyset}  {\mathcal S}_j,
\end{equation}
respectively.
Note that definition \cref{eqn:region} includes all of the elements $\Tm_j$ that share either a face or a vertex with element $\Tm_i$ -- we refer to all of the elements in ${\mathcal R}_i$ as 
{\it vertex-neighbors} of element $\Tm_i$; a 2D version is depicted in \cref{fig:RIDG_2D} .
Also note that on a Cartesian mesh in $\mdim$ dimensions, region ${\mathcal R}_i$ contains exactly $3^{\mdim}$ spacetime elements. On each spacetime element we write the predicted solution as
\begin{equation}
\label{eqn:pred_ansatz}
\vec{q}\left(t,\vec{x}\right) \Bigl|_{{\mathcal S}_i} \approx w_{i} := 
\vec{\Psi}^T \mat{W}_i,
\end{equation}
where $\mat{W}_i \in \reals^{\mpred \times \meq}$, $\vec{\Psi} \in \reals^{\mpred}$, 
\begin{equation*}
\mpred := \left(\mdeg+1 \right)^{\mdim+1},
\end{equation*}
 and $\vec{\Psi}$ is an orthonormal basis on 
the spacetime element:
\begin{equation*}
\frac{1}{\bigl|{\mathcal S}_i \bigr|} \int_{S_i} \vec{\Psi} \, \vec{\Psi}^T \, d\tau \, d\xi  = \mat{\mathbb I} \in \reals^{\mpred\times\mpred}.
\end{equation*}

In order to arrive at an algebraic system of equations to solve for the unknown spacetime coefficients in region ${\mathcal R}_i$, we multiply \cref{eqn:conslaw}
by test functions $\vec{\Psi}$ and then integrate over each ${\mathcal S}_j \in
{\mathcal R}_i$. We apply integration-by-parts in both space and time. The integration-by-parts in time connects the current time-slab $[t^n, t^{n+1}]$ to the solution in the previous time-slab (i.e., this enforces causality). In space we treat element boundaries differently depending on whether or not those boundaries are (1) strictly internal to region ${\mathcal R}_i$ or (2) on the boundary
of region ${\mathcal R}_i$. For strictly internal boundaries we evaluate the resulting surface integral through standard upwind fluxes (i.e., Rusanov fluxes) that utilize the solution on both sides of the interface. However, on
the boundary  of ${\mathcal R}_i$ we evaluate the resulting surface integrals
only using the trace of the solution that is internal to ${\mathcal R}_i$.
The effect of these choices is to insulate region ${\mathcal R}_i$ from all elements exterior to ${\mathcal R}_i$. Solving the algebraic equations in
region ${\mathcal R}_i$ resulting from the above described integration-by-parts,
gives us the {\it regionally-implicit prediction}.
This setup is depicted in the 2D setting in \cref{fig:RIDG_2D} .

Algebraically, the system that must be solved in each region has the following form
for each $j: \Tm_j \in {\mathcal R}_i$ and for each $k = 1,2, \ldots, \mpred$:
\begin{align}
\label{eqn:pred_full_system}
\begin{split}
\vec{R}_{jk} :=& \underbrace{\left[ \int_{{\mathcal T}_j}   \left( \Psi_k
\vec{\Psi}^T \right)\Bigl|_{t^{n+1}}   \, d\vec{x} -\int_{t^n}^{t^{n+1}} \int_{\Tm_j}  \Psi_{k,t} \, \vec{\Psi}^T \,  d\vec{x} \, dt \right] \, \mat{W}_j}_{\text{(time term)}} \\
-& \underbrace{\int_{t^n}^{t^{n+1}} \int_{\Tm_j} 
 \vec{\nabla} \Psi_k \cdot \mat{F}\left( \vec{\Psi}^T \, \mat{W}_j \right) \, d\vec{x} \, dt}_{\text{(internal flux term)}} \\ 
+& \underbrace{\int_{t^n}^{t^{n+1}} \oint_{\partial \Tm^{\star}_j}  \Psi_k \, 
\vec{{\mathcal F}}\left( \vec{w}^h_{+}, \vec{w}^h_{-};
\vec{n} \right) \, d\vec{s} \, dt}_{\text{(surface flux term)}} - 
\underbrace{\left[ \int_{{\mathcal T}_j} \left( \Psi_k
\vec{\Phi}^T \right)\Bigl|_{t^{n}} \, d\vec{x} \right] \, \mat{Q}^{n}_j}_{\text{(causal source)}} = \vec{0},
\end{split}
\end{align}
where $\partial \Tm^{\star}_j$ is the part of the boundary of $\Tm_j$ that
is interior to ${\mathcal R}_i$ and $\vec{n}$ is an outward pointing normal
to on $\partial \Tm^{\star}_j$. The states $\vec{w}^h_{\pm}$ are the solution on
either sides of the boundary $\partial \Tm^{\star}_j$. In this work, the numerical flux, $\vec{{\mathcal F}}$, is taken to be the Rusanov (sometimes called the local Lax-Friedrichs) flux \cite{article:Ru61}. 

On a Cartesian mesh in $\mdim$ dimensions, equation \cref{eqn:pred_full_system} represents a nonlinear algebraic system of size $\meq \cdot \mpred \cdot 3^\mdim$. In fact, the only portion of the solution to \cref{eqn:pred_full_system} that we actually need to retain is the solution on the central element: $\Tm_i$; the solution on the remaining elements in ${\mathcal R}_i$ will be discarded (again, see \cref{fig:RIDG_2D} ). 

This prediction step portion of the RIDG scheme is by far the most expensive. Furthermore, the solution to \cref{eqn:pred_full_system}  by itself is quite useless: this solution is not even consistent with the original PDE due to the fact that region ${\mathcal R}_i$ has been insulated from all other regions in the computational domain. However, using this solution as a prediction inside of an appropriate correction step produces a scheme that is both high-order accurate and stable up to CFL numbers that significantly exceed Lax-Wendroff DG and SSP-RKDG \cite{RIDG_paper_2019}.

\subsection{Correction step for RIDG}
For the correction step we assume a solution of the form \cref{eqn:phi_basis_1d}
and use of a version of equation \cref{eqn:semi_discrete_dg} where the solution
inside each of the spacetime integrals is replaced by the predicted solution
\cref{eqn:pred_ansatz} :
\begin{equation*}
\begin{split}
\vec{Q}_{i}^{k}\left( t^{n+1} \right) &= 
\vec{Q}_{i}^{k}\left( t^{n} \right)
+ \frac{1}{|\Tm_i|} \int_{t^n}^{t^{n+1}} \int_{\Tm_i} 
\grad \Phi_{k} \cdot \mat{F}\left(
 \vec{\Psi}^T \, \mat{W}_i \right)  \, d\vec{x} \, dt \\
&- \frac{1}{|\Tm_i|}\int_{t^n}^{t^{n+1}} \oint_{\partial \Tm_i} \Phi_k \, \vec{{\mathcal F}}\left(\vec{w}^h_{+}, \vec{w}^h_{-};
\vec{n} \right) \, d\vec{s} \, dt.
\end{split}
\end{equation*}
In the above equatoion $\vec{n}$ is an outward-pointing normal vector to $\partial \Tm_i$, $\vec{w}^h_{+}$ and $\vec{w}^h_{-}$ are the states on either side of the boundary $\partial \Tm_i$, and $\vec{{\mathcal F}}$ is the numerical flux. In this work we always make use of the Rusanov (sometimes called the local Lax-Friedrichs) flux
\cite{article:Ru61}. 
This portion of the regionally-implicit DG update is fully explicit and computationally inexpensive relative to the prediction step. 

%% file: tikz/ridg2d.tex
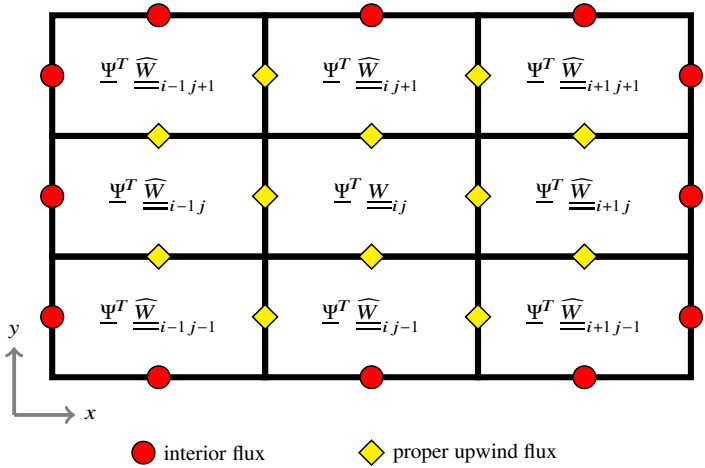
\begin{figure}[t]
\centering
\begin{tikzpicture}[scale=0.4, 
circ/.style={scale=1.5, shape=circle, inner sep=2pt, draw, fill=red,node contents=},
dia/.style={scale=1.2, shape=diamond, inner sep=2pt, draw, fill=yellow,node contents=},
amp/.style = {scale=0.75, regular polygon, regular polygon sides=3, shape border rotate=-180,
              draw, fill=yellow, inner sep=2pt, node contents=}]

\draw[->,line width=0.5mm, gray] (1.750,-1.25)--(1.750,1) node[above,black]{$y$};
\draw[->,line width=0.5mm, gray] (1.750,-1.25)--(3.75,-1.25) node[right,black]{$x$};

\draw [line width=0.8mm] (3,0) -- (10,0) -- (10,4) -- (3,4) -- cycle;
\draw [line width=0.8mm] (10,0) -- (17,0) -- (17,4) -- (10,4) -- cycle;
\draw [line width=0.8mm] (17,0) -- (24,0) -- (24,4) -- (17,4) -- cycle;

\draw [line width=0.8mm] (3,0+4) -- (10,0+4) -- (10,4+4) -- (3,4+4) -- cycle;
\draw [line width=0.8mm] (10,0+4) -- (17,0+4) -- (17,4+4) -- (10,4+4) -- cycle;
\draw [line width=0.8mm] (17,0+4) -- (24,0+4) -- (24,4+4) -- (17,4+4) -- cycle;

\draw [line width=0.8mm] (3,0+8) -- (10,0+8) -- (10,4+8) -- (3,4+8) -- cycle;
\draw [line width=0.8mm] (10,0+8) -- (17,0+8) -- (17,4+8) -- (10,4+8) -- cycle;
\draw [line width=0.8mm] (17,0+8) -- (24,0+8) -- (24,4+8) -- (17,4+8) -- cycle;

\draw [line width=0.2mm, black] node (r2) at (3, 2) [circ];
\draw [line width=0.2mm, black] node (r5) at (10, 2) [dia];
\draw [line width=0.2mm, black] node (r8) at (17, 2) [dia];
\draw [line width=0.2mm, black] node (r11) at (24, 2) [circ];
\draw [line width=0.2mm, black] node (r14) at (6.5, 4) [dia];
\draw [line width=0.2mm, black] node (r17) at (6.5, 0) [circ];
\draw [line width=0.2mm, black] node (r20) at (13.5, 4) [dia];
\draw [line width=0.2mm, black] node (r23) at (13.5, 0) [circ];
\draw [line width=0.2mm, black] node (r26) at (20.5, 4) [dia];
\draw [line width=0.2mm, black] node (r29) at (20.5, 0) [circ];
\draw [line width=0.2mm, black] node (r29) at (20.5, 12) [circ];
\draw [line width=0.2mm, black] node (r29) at (13.5, 12) [circ];
\draw [line width=0.2mm, black] node (r29) at (6.5, 12) [circ];
\draw [line width=0.2mm, black] node (r29) at (20.5, 8) [dia];
\draw [line width=0.2mm, black] node (r29) at (13.5, 8) [dia];
\draw [line width=0.2mm, black] node (r29) at (6.5, 8) [dia];
\draw [line width=0.2mm, black] node (r29) at (10, 6) [dia];
\draw [line width=0.2mm, black] node (r29) at (10, 6+4) [dia];
\draw [line width=0.2mm, black] node (r29) at (17, 6) [dia];
\draw [line width=0.2mm, black] node (r29) at (17, 6+4) [dia];
\draw [line width=0.2mm, black] node (r29) at (3, 6+4) [circ];
\draw [line width=0.2mm, black] node (r29) at (3, 6) [circ];
\draw [line width=0.2mm, black] node (r29) at (24, 6+4) [circ];
\draw [line width=0.2mm, black] node (r29) at (24, 6) [circ];

\node[] at (6.5,2) {$\vec{\Psi}^T \, \mat{\widehat{W}}_{\, i-1 \, j-1}$};
\node[] at (13.5,2) {$\vec{\Psi}^T \, \mat{\widehat{W}}_{\, i \, j-1}$};
\node[] at (20.5,2) {$\vec{\Psi}^T \, \mat{\widehat{W}}_{\, i+1 \, j-1}$};

\node[] at (6.5,2+4) {$\vec{\Psi}^T \, \mat{\widehat{W}}_{\, i-1 \, j}$};
\node[] at (13.5,2+4) {$\vec{\Psi}^T \, \mat{W}_{ij}$};
\node[] at (20.5,2+4) {$\vec{\Psi}^T \, \mat{\widehat{W}}_{\, i+1 \, j}$};

\node[] at (6.5,2+8) {$\vec{\Psi}^T \, \mat{\widehat{W}}_{\, i-1 \, j+1}$};
\node[] at (13.5,2+8) {$\vec{\Psi}^T \, \mat{\widehat{W}}_{\, i \, j+1}$};
\node[] at (20.5,2+8) {$\vec{\Psi}^T \, \mat{\widehat{W}}_{\, i+1 \, j+1}$};

\draw [line width=0.2mm, black] node (l2) at (6, -2.5) [circ,label=right:{interior flux}];
\draw [line width=0.2mm, black] node (l3) at (13.5, -2.5) [dia,label=right:{proper upwind flux}];

\end{tikzpicture}
\caption{Shown is the 2D Cartesian region ${\mathcal R}_i$ over which the
RIDG prediction step is carried out. In the RIDG prediction step, all of the states, excepting only the one belonging to the middle element, are only temporary variables
and will be discarded once the predicted solution in element $ij$ has been computed -- to make note of this we place hats over the temporary variables.\label{fig:RIDG_2D}}
\end{figure}

%% file: secJacobians.tex

\section{Efficient implementation of RIDG via quasi-quadrature-free Jacobian matrix assembly}
\label{sec:jacobians}
As described in \cref{sec:RIDG1D_predict} , the prediction step in the regionally-implicit DG method is the most computationally expensive, since it requires the solution of a nonlinear system of algebraic equations over each region ${\mathcal R}_i$. Due to the computational complexity of space-time DG methods, efficient matrix assembly is required if high order RIDG methods are to be viable. Therefore, in order to make this step as efficient as possible, we develop in this section
a quasi-quadrature-free Jacobian matrix assembly. 
 
Quadrature-free schemes have provided efficient alternatives to quadrature in DG schemes for many years \cite{Atkins1998}.  Depending on the basis choice, they can benefit from reduced computational complexity for high-dimensional problems without sacrificing accuracy.  That is, the quadrature free schemes can be designed to inherit the same assumptions about an underlying integrand as do discrete quadrature schemes. Thus, quadrature-free integrals are efficient replacements for standard quadrature schemes, especially when orthogonal bases are used.  Matrix assembly for nodal elements can have additional exploits \cite{Engsig-Karup2016}.  With some work, quarature-free schemes can also be accurately extended to unstructured geometries \cite{chan2016}. 
In the current work we extend this previous work to high-order space-time Jacobian matrix assembly, as required by the Newton method for solving the nonlinear systems inherent to the RIDG prediction step \cite{Guthrey2017,RIDG_paper_2019}.   

The prediction step algebraic system given by \cref{eqn:pred_full_system} is typically solved via a Newton's method, which can be written as
\begin{equation*}
 \vec W^{(m+1)} = \vec W^{(m)} - \left[{\mat{J}}\left(\vec W^{(m)}\right)\right]^{-1} \vec{R}
 \left(\vec W^{(m)}\right),
\end{equation*}
where $\vec{R}$ is shorthand for the residual defined in \cref{eqn:pred_full_system} written as a column vector, $\vec{W}$ is shorthand for all of the unknown cofficients
in region ${\mathcal R}_i$ written as a column vector, $m$ is the Newton iteration counter, and $\mat{J}$ is the Jacobian of $\vec{R}$ with respect to $\vec{W}$.
   
   For ease of discussion, let us focus our description of the Jacobian-free
   implementation on only the {\it internal flux term} in the residual and let us consider only a scalar case ($\meq=1$) -- these restrictions can easily
   be removed:
\begin{equation}\label{residual}
\begin{split}
     {R}_{jk} &= (\text{time term}) -  \int_{t^n}^{t^{n+1}} \int_{\Tm_j} 
 \vec{\nabla} \Psi_k \cdot \vec{F}\left( \vec{\Psi}^T \, \vec{W}_j \right) \, d\vec{x} \, dt \\ &+  (\text{surface flux term}) - (\text{causal source}),
 \end{split}
\end{equation}
for each $j: \Tm_j \in {\mathcal R}_i$ and for each $k = 1,2, \ldots, \mpred$, 
where $\vec{W}_j \in \reals^{\mpred}$ are the unknown coefficients in spacetime element  ${\mathcal S}_j$.
The Jacobian of this residual with respect to the degrees of freedom $\vec W$ is 
\begin{equation}
\label{Jacobian}
\begin{split}
     {J}_{jk\ell} &= \frac{\partial}{\partial {W_{\ell}}}(\text{time term}) -  \int_{t^n}^{t^{n+1}} \int_{\Tm_j} 
 \vec{\nabla} \Psi_k \cdot \vec{F}'\left( \vec{\Psi}^T \, \vec{W}_j \right) \, {\Psi}_{\ell} \, d\vec{x} \, dt
  \\ &+  \frac{\partial}{\partial {W_{\ell}}}(\text{surface flux term}),
     \end{split}
\end{equation}
where $\vec{F}'(q) \cdot \vec{n}$ is the flux Jacobian of the hyperbolic conservation law in direction $\vec{n}$. 

In order to assemble the above Jacobian matrix for the Newton iteration, several methods are available. 
\begin{itemize}
\item {\bf Quadrature:} We could perform the integrals in \cref{Jacobian} via quadrature for every $k,\ell \leq \theta$, where $\theta$ is the number of terms in our space-time Legendre basis.  This method is able to take advantage of the sparsity pattern of our block-stencil.  The computational complexity of such a method is
    \begin{equation} \label{quadcomplexity} 
    \text{timing} \approx \mathcal O( \theta^2 (\mdeg+1)^{\mdim+1} ) \approx 
     \mathcal O( (\mdeg+1)^{3\mdim+3} ) , 
    \end{equation}
    which will become incredibly expensive for high-order methods in three dimensions.

\item {\bf Perturbation:} One alternative to direct quadrature is to
 approximate the Jacobian ${J}_{jk\ell}$ via either finite differences or Gateaux derivatives of the residual \cref{residual} ; this is also known as a perturbation method.  This method is very simple to implement, as you only need to define the residual for your Newton iteration.  However, since the residual needs to be recomputed a number of times equal to the number of degrees of freedom, the computational complexity is roughly the same as \cref{quadcomplexity} .

\item {\bf Quadrature-free:} If the flux function is sufficiently simple
(i.e., a polynomial function of the solution), then an alternative to the quadrature and perturbation
approaches is to perform exact integrations of the terms needed in the Jacobian; this was the method employed in  Guthrey and Rossmanith \cite{RIDG_paper_2019}.  
Unfortunately, this method does not generalize well, even to relatively simple rational flux functions such as those seen in the compressible Euler equations.  

\end{itemize}

In this work we consider an alternative to all of these approaches; it is similar to the quadrature-free approach, but applies to general fluxes (i.e., non-polynomial).  The idea is this: we project components of the flux Jacobian, $\vec{F}'(w)$, onto the  polynomial basis $\vec{\Psi}$:
\begin{equation}\label{Jacproj}
  \vec F'_{\, p} =  \left\{ \int_{t^n}^{t^{n+1}} \int_{\Tm_j}  \Psi_p \, \vec F'(w) \, d\vec x \, dt \right\} \biggl/ \left\{ \int_{t^n}^{t^{n+1}} \int_{\Tm_j}  \left( \Psi_p \right)^2  \, d\vec x \, dt \right\}.
\end{equation}
Using the expansion of $\vec F'(w)$ in the expression for the residual \cref{Jacobian} , we obtain
\begin{equation}\label{Jacobian2}
\begin{split}
     {J}_{jk\ell} &= \frac{\partial}{\partial {W_{\ell}}}(\text{time term}) -  \sum_{p=1}^{\infty}\left( \int_{t^n}^{t^{n+1}} \int_{\Tm_j} 
 \vec{\nabla} \Psi_k \, \Psi_p \, {\Psi}_{\ell} \, d\vec{x} \, dt \right)
 \cdot \vec{F}'_{\, p}
  \\ &+  \frac{\partial}{\partial {W_{\ell}}}(\text{surface flux term}).
     \end{split}
\end{equation}
Due to the orthogonality of our basis, there exists some $L$ such that for all $p > L$ and all needed $j$, $k$, and $\ell$, the integrals in the above expression all vanish (see Gupta and Narasimhan \cite{Gupta2007}).
That is, there are a limited amount of the integrals in \cref{Jacobian2} that are actually nonzero.  
Furthermore, these integral expressions can be precomputed exactly and stored, and quadrature is only needed for the projection \cref{Jacproj} .  
Thus, for each $j,k,\ell$ we may form a list $\alpha_{jk\ell}$ of the indices $p$ such that the associated integral was nonzero and a list $\vec\beta_{j k \ell p}$ of the precomputed results of the nonzero integrals. Finally, we can then simplify \cref{Jacobian2} to 
the following:
\begin{equation*}
\begin{split}
     {J}_{jk\ell} &= \frac{\partial}{\partial {W_{\ell}}}(\text{time term}) -  \sum_{p=1}^{\alpha_{jk\ell}}\left( \vec\beta_{j k \ell p}
 \cdot \vec{F}'_{\, p} \right)
  +  \frac{\partial}{\partial {W_{\ell}}}(\text{surface flux term}).
     \end{split}
\end{equation*}
\input{tables/qqf}
This quasi-quadrature-free method maintains the same order of accuracy as the
 quadrature or perturbation strategies described above, but in practice is orders of magnitude faster, because the space-time quadrature is performed once per degree of freedom per element as opposed to once per square of the degrees of freedom per element.  We compare timings of these three methods with
 $\mdim=3$ in Table \cref{table:qqf} .  At this time we do not have an {\it a priori} estimate for the computational complexity of the quasi-quadrature-free method, but we are able to compute {\it a posteriori} estimates for each method: 
\begin{align*}
    \text{perturbation timings} \approx  & \; \mathcal O( (\mdeg+1)^{11.9} )  \\
    \text{quadrature timings} \approx  & \; \mathcal O( (\mdeg+1)^{11.6})  \\
    \text{quasi-quadrature-free timings} \approx  & \;\mathcal O( (\mdeg+1)^{9.3} ) 
\end{align*}
We see that the perturbation and traditional quadrature routines exhibit the expected computational complexity described in \cref{quadcomplexity}, but that the quasi-quadrature-free method exhibits a significantly lower computational complexity.

For brevity and clarity we have omitted the discussion of the temporal derivative terms and surface flux terms. However, we note that the time terms of the matrix assembly of $J_{jk\ell}$ can be precomputed and stored, and no projection is needed during the matrix assembly.  Furthermore, the surface terms in general involve a numerical flux, for which a numerical flux Jacobian must be generated or approximated before it is projected.  These types of terms affect several blocks of the Jacobian matrix, exactly in accordance to the block-stencil derived from our numerical flux.

To review, we have efficiently extended quadrature-free schemes to our space-time prediction procedure needed for our regionally-implicit DG method. This implementation offers incredible speedups to a very computationally expensive method in a way that is extensible to arbitrary problems and relatively simple to implement. We will use this strategy in the remainder of this paper.

%% file: tables/qqf.tex
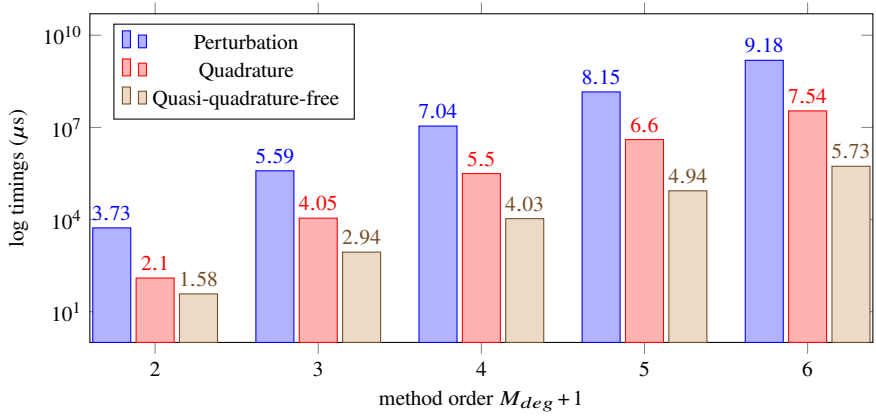
\begin{figure}
    \centering
\pgfplotstableread[row sep=\\,col sep=&]{
    interval & pert & quad & qqf \\
2	&	5357	&	125	&	38	\\
3	&	386006	&	11100	&	871	\\
4	&	11059200	&	314000	&	10700	\\
5	&	142762500	&	4010000	&	86300	\\
6	&	1529150400	&	34400000	&	543000	\\
    }\mydata
    \begin{tikzpicture}
    \begin{axis}[
            ymode = log,
            log basis y={10},
            ybar,
            bar width=.5cm,
            width=\textwidth,
            height=.5\textwidth,
            legend pos=north west,
            symbolic x coords={2,3,4,5,6},
            xtick=data,
            nodes near coords,
            nodes near coords align={vertical},
            ymin=1,ymax=5e10,
            ylabel={log timings ($\mu$s)},
            xlabel={method order $M_{deg}+1$}
            ]
        \addplot table[x=interval,y=pert]{\mydata};
        \addplot table[x=interval,y=quad]{\mydata};
        \addplot table[x=interval,y=qqf]{\mydata};
        \legend{Perturbation,Quadrature,Quasi-quadrature-free} \end{axis}
\end{tikzpicture}
    \caption{Newton iteration Jacobian matrix assembly time for one 3D+1 space-time prediction for a nonlinear problem, using (1, blue) a theoretical estimate of the cost of a perturbation method obtained by multiplying the timings of the residual computation by the number of degrees of freedom in the Newton iteration (2, red) traditional quadrature routines to compute the integrals in \cref{Jacobian}, and (3, tan) the quasi-quadrature-free routine discussed in this appendix. Runtimes are $\log_{10}$ of averages of 10 trial runs in microseconds.  We see that the quasi-quadrature-free method is the most efficient, and scales better as we increase $\mdeg$.  For $\mdeg+1=6$, the quasi-quadrature-free method is almost two orders of magnitude faster than the traditional quadrature scheme.
    }
    \label{table:qqf}
\end{figure}

%% file: secterminology.tex

\section{Basic terminology}
\label{sec:terminology}
The goal of this work is to extend the convergence studies performed in Guthrey and Rossmanith \cite{RIDG_paper_2019} to the case of strong scaling in an HPC setting. 
Before reporting results in subsequent sections, we describe here some of the key terms needed to analyze the efficiency and accuracy of parallelized DG schemes.
We also note that here and throughout the remainder of the paper we consider only the scalar: $\meq=1$. All the findings in this work generalize to the more general case, but for brevity and clarity we assume for the remainder of this current work that $\meq=1$.

We define the following key terms.
\begin{description}
\item[{\bf mesh}:] The description of the number of elements.  All methods considered in this paper use a uniform Cartesian mesh.

\item {\bf dof} : The total number of degrees of freedom. We compute this by multiplying basis element size ($\theta$) by the number of  elements ({\bf mesh}):
\begin{equation*}
\text{\bf dof}( \theta ,\text{\bf mesh} ) = \theta\times\text{\bf mesh}.
\end{equation*}    

\item[{\bf efom}:] High-order DG methods can be compared to first-order finite volume methods by (among other ways) comparing the equivalent number of the total degrees of freedom. Thus, {\bf efom},
equivalent first order mesh, indicates what size of a Cartesian mesh with one degree of freedom per cell would have the equivalent number of degrees of freedom for a given mesh and method order. The {\bf efom} for a method with a basis of size $\theta$ for a $\mdim$-dimensional problem is computed using 
\begin{equation}
\label{eqn:EFOM}
\text{\bf efom}( \theta ,\text{\bf mesh} ) = \left[ \paren{\text{\bf dof}}^{\frac{1}{\mdim}}\right]^{\mdim}
=  \left[\paren{\theta\times\text{\bf mesh}}^{\frac{1}{\mdim}}\right]^{\mdim}.
\end{equation}
For example, a $10\times 10$ mesh in 2D with a basis of size $16$ has an {\bf efom} of $40 \times 40$ since it has the equivalent number of degrees of freedom as a $40\times 40$ first order mesh. 

\item [{\bf error}:] The approximate relative error in the $L^2$ norm of our numerical solution.  This error is obtain by considering an exact projection $q_{\text{exact}}$ onto an infinite basis, $\vect \Phi$, on each element: 
    $\vec{X}_i := \proj{\vect \Phi^\infty} q_{\text{exact}}$.  Comparing
    this exact expansion to the numerical solution gives the
    following element-wise error definition:
    \begin{equation*}
     e(\vec{x})\Bigl|_{\Tm_i} :=  \left| \underbrace{\sum\limits_{\ell = 1}^\infty  \Phi_\ell(\vec{x}) \, X^{\ell}_{i}}_{\text{exact}} - \underbrace{\sum\limits_{\ell = 1}^{\mbasis}  \Phi_{\ell}(\vec{x}) \, Q^{\ell}_i }_{\text{numerical}}\right|.
    \end{equation*}
    Using the orthonormality of the basis functions
     gives the following $L^2$ error:
        \begin{equation*}
     \norm[L^2]{e}^2 =  \sum_{i} \left [\underbrace{\sum\limits_{\ell = 1}^{\mbasis}  \left( X_{i}^{\ell} - Q^{\ell}_i  \right)^2}_{\text{error in coeffs}}
    +  \underbrace{\sum\limits_{\ell = \mbasis+1}^{\mbasis^{+}} \left(X_{i}^{\ell} \right)^2}_{\text{dominant trunc. error}} + \underbrace{\sum\limits_{\ell = \mbasis^{+}+1}^\infty \left( X_{i}^{\ell}\right)^2 }_{\text{high-order trunc. error}}\right],
    \end{equation*}
    where $\mbasis$ and $\mbasis^{+}$ are the number of basis functions in 
    ${\mathbb Q} \left(\mdeg, \mdim \right)$ and 
    ${\mathbb Q} \left(\mdeg+1, \mdim \right)$, respectively.
    To sufficient accuracy the error can be approximated by 
    discarding all the terms past  $\ell=\mbasis^{+}$, resulting
    in the following approximate relative error:
    \begin{equation}
    \label{eqn:l2relerror}
    \text{ \bf error}  = \frac{\norm[L^2]{e}}{\norm[L^2]{q_{\text{exact}}}}  \approx \sqrt{ \frac{
    \sum\limits_{i=1} \left[ \sum\limits_{\ell = 1}^{\mbasis}  \left( X_{i}^{\ell} - Q^{\ell}_i  \right)^2 + \sum\limits_{\ell=\mbasis+1}^{\mbasis^{+}} \left(X_i^{\ell} \right)^2 \right]
    }{
    \sum\limits_{i=1} \left[ \sum\limits_{\ell = 1}^{\mbasis^{+}} \left( X_i^{\ell} \right)^2 \right]
    }}.
    \end{equation}
    \item {\bf aproximate order of accuracy}: For the convergence studies performed in this paper, we approximate the convergence rate $M$ using
    \begin{equation}
\label{eqn:Mratio}
\text{\bf error}(h) = c h^M + {\mathcal O}\left(h^{M+1} \right) \quad \Longrightarrow \quad
M \approx \frac{\log\paren{{\text{\bf error}(h_1)}/{\text{\bf error}(h_2)}}}{\log\paren{{h_1}/{h_2}} },
\end{equation}
where $h$ is the mesh spacing.

\item[{\bf runtime}:]  For each experiment we compute the wall clock runtime of each method in seconds. For RIDG methods, the computational effort is dominated by small dense matrix inverses in the prediction step, and thus for a method that requires $N_t$ time-steps with a mesh of size $N^\mdim$ and a space-time basis of size $\theta_T$, we expect the runtime to scale as
\begin{equation}
 \text{\bf runtime} \approx \mathcal O\left( N_t \cdot N^\mdim \cdot (3^\mdim \theta_T)^3\right).
 \end{equation}
If we use the $\mathcal Q$ spacetime basis for the prediction step as discussed in \cite{RIDG_paper_2019}  and we let $M = \mdeg+1$, then $\theta_T = M^{\mdim+1}$ and so for explicit time-stepping where $\Delta t \approx  \nu \Delta x \approx \nu/N$, where $\nu$ is the CFL number, we get that   
\begin{equation}
\label{eqn:scaleruntime}
 \text{\bf runtime}\approx \mathcal O\left( 27^{\mdim} \cdot \nu^\inv \cdot
   N^{\mdim+1} \cdot M^{3\mdim+3}  \right).
 \end{equation}
RIDG methods have the fortunate property that $\nu^\inv \approx \mathcal O(1)$, whereas methods such as SSP-RKDG experience  $\nu^\inv \approx \mathcal O(M^2)$.

\item[{\bf quality}:] When comparing various methods in terms of efficiency, the two most prevalent metrics are {\bf runtime} and {\bf error}.  One could consider how fast various methods reach a fixed error (time to solution), or vice versa. In this paper we combine these two metrics into a third metric we call {\bf quality}.  We define the quality of the solution as
\begin{equation}
 \text{\bf quality} = -\log\paren{ \text{\bf error}\times\text{\bf runtime}}.
\label{eqn:Quality}
\end{equation}
We note that this metric is on a logarithmic scale, so changes of $\pm 1$ for this metric are quite significant.  If we assume that RIDG methods of order $M = \mdeg+1$, have an error convergence scaling of $\mathcal O(N^{-M})$, then using \cref{eqn:scaleruntime} we expect the quality to scale as 
\begin{equation}
 \text{\bf quality} \approx \mathcal O\left( \nu + (M-\mdim-1)\log(N) - (3\mdim+3)\log(M) \right).
\label{eqn:ridgqualityraw}
\end{equation}
This provides the intuition that high order RIDG methods will have a high quality, offset by increased runtime costs associated with the dense matrix inverses. However, as we consider higher mesh resolutions, the quality of a method of fixed order will be driven higher by the shrinking error.  For smaller mesh resolutions, the quality will be adversely dominated by the high runtime cost of the methods.

\item [{\bf dof/c}:] This is the number of degrees of freedom per compute core, and is computed as:
\begin{equation}
    \label{eqn:ndofspercore}
\text{\bf dof/c} =  
\frac{ \theta\times {\bf mesh}  }{ \bf cores}
= \frac{ \theta\times {\bf mesh}  }{ \bf tasks \times cores/task} ,
\end{equation}
where $\theta$ is the size of the DG basis on each element. This quantity is one way to measure the relative compute intensity per core for a given experiment. 

\item [{\bf speedup} over single task runtimes:] We measure the ratios of runtimes of the single task run versus the multi-task runs to compare how much faster the latter case runs. Ideally we expect that the speedup is equal to the total number of tasks being used. It can be described with the formula 
\begin{equation}
    \label{eqn:Speedup}
\text{\bf strong speedup}( \text{\bf tasks} = i )  =  \frac{ \text{\bf runtime}( \text{\bf tasks} = 1) }{ \text{\bf runtime}( \text{\bf tasks} = i)  }.
\end{equation}
Note that this definition of speedup compares runtimes for a given set of tasks as opposed to cores.  That is, we do not compare runtimes against serial execution. 

\item[{\bf strong efficiency}:] This is a measure of how our speedup is to the ideal case. It is simply the ratio of actual speedup versus ideal speedup. As we add tasks, the parallel efficiency for a strong scaling test is:
\begin{equation}
    \label{eqn:Efficiency}
\text{\bf strong efficiency}( \text{\bf tasks} = i )  =  \frac{ \text{\bf speedup}( \text{\bf tasks} = i) - 1 }{  i - 1  }\times 100\% .
\end{equation}

An efficiency of $100\%$ corresponds to perfect linear speedup scaling with the number of tasks (or cores). An efficiency of $0\%$ means that as we used more tasks, the runtime remained unchanged or increased (we map negative efficiencies to $0\%$).  

\item [{\bf comms}:] This is an estimate for the total number of MPI communications. It is estimated using the formula 
\begin{equation}
\label{eqn:comms}
\text{\bf comms}  = {\bf tasks}\times {\bf timesteps}\times {\bf stages/timestep}\times {\bf comms/stage} .
\end{equation}
We approximate the number of timesteps via the CFL relation 
\begin{equation}
 \Delta t = \nu \Delta x \quad \implies \quad  \text{\bf timesteps} = \nu^{-1} N,
 \end{equation}
where $N$ is the number of mesh elements in each direction. We note that RKDG methods of $M_{deg}=3$ have 9 stages per timestep, and that all RIDG methods have 2 stages (predictor and corrector).  There are 8 communications per stage in 2D and 26 in 3D.

\end{description}

%% file: sec_implement.tex

\section{Parallel implementation of RIDG}
\label{sec:parallel}


We now briefly discuss the domain decomposition strategy implemented for achieving high performance computing with RKDG and RIDG methods.  We consider a Cartesian mesh evenly subdivided into a $N\times N$ (2D) or $N\times N\times N$ (3D) grid of Cartesian  submeshes, where $N^2$ (2D) or $N^3$ (3D) is the number of compute nodes used. Each node (i.e., submesh) is associated with a single MPI task. 
OpenMP is used for shared-memory parallelization of operations on each Cartesian submesh. MPI is used to perform data communications across the Cartesian submesh interfaces, often called the {\it ghost zone} or {\it halo region} communications.  Submeshes that share interfaces will need to communicate multiple times per timestep, as detailed in the next two sections.

\subsection{Domain decomposition strategy for RKDG}
In order to compute the needed fluxes for each stage of RKDG, each Cartesian submesh requires information from face-neighbor submeshes. This means that for RKDG the submeshes must communicate with neighboring submeshes in each stage.
The communication latency for this operation can be hidden by using non-blocking MPI routines for the intra-face communications before we perform the volume integrals over the submesh. Each volume integral has a relatively high arithmetic intensity and can be completed using information completely local to each element.  Once these volume integrals are computed, we simply wait until intra-face communications are completed, which ideally has already occurred. Then, the fluxes and thus the update for the RK stage can be computed. 
Lastly, each timestep requires an all-to-all communication to enforce the CFL time-step restriction by communicating the CFL number used by each submesh.

\subsection{Domain decomposition strategy for RIDG}
For RIDG the prediction step requires that each given Cartesian submesh has information from vertex-neighbor submeshes, since forming a region for a cell involve the cell's vertex-neighbors.  Latency associated with this communication can be hidden by using non-blocking MPI routines for the interface communications. After these communications are started, we perform the prediction step for elements not on the submesh boundary, as data for their vertex neighbors is located in shared memory.  
Once the boundary communications are complete, we can and then compute the predictions for boundary elements. We must communicate these boundary predictions as they are needed for the correction step. Again we may use non-blocking routines and finish computing the non-boundary predictions during these communications. Once all predictions are formed and the communications are complete, we continue to the correction step. Just as with RKDG, each timestep requires an all-to-all communication to enforce the CFL time-step restriction by communicating the CFL number used by each submesh. This all-to-all communication can also be performed using a non-blocking MPI routine while the correction update is computed.  For nonlinear problems, the prediction step and its associated communications must be repeated until the region residual is driven down to some tolerance. 


%% file: section_1d.tex

\section{Numerical results}
In order validate the proposed implementation of the regionally-implicit scheme we consider in this section examples 1D, 2D, and 3D. The 1D code runs sufficiently fast that no parallelization is required; we provide the 1D results mostly to 
demonstrate the efficiency gains of high-order RIDG.
Subsequently we consider 2D and 3D examples, all of which are implemented in parallel.

\subsection{1D results}
\label{sec:1dresults}
We consider the following 1D linear advection equation:
\begin{equation}
\label{eqn:adv1d}
\piecewise{    q_{,t} + q_{,x}  = 0 & \text{for} \quad (t,x) \in  \left [ 0, T \right] \times \left [0,1 \right], \\ 
q(t=0,x) = \exp{((x-c)^2-\omega^2 )^{-1}} & \text{if} \quad (x-c)^2 < \omega^2, \\
q(t=0,x) = 0  & \text{if} \quad (x-c)^2 \geq \omega^2,} 
\end{equation}
with $c=\frac{1}{2}$ and $\omega = \frac{1}{3}$. Note that the initial condition is clearly $C^\infty$ for $\abs{x-c}<\omega$; furthermore, it can be shown that all derivatives for this function vanish as $(x-c)^2 \to \omega^2$ from the left, and hence the initial condition is actually $C^\infty$ over all of $\reals$.  This level of smoothness is required to test the convergence properties of arbitrarily high order methods. 

\input{tables/advection1d.tex}

The results of these convergence tests are shown in  \cref{table:1D} . 
We notice that the maximum allowable CFL restriction for the RIDG methods is bounded from below by some $\nu_{\text{min}}$; the key point is that $\nu_{\text{min}}$ for RIDG is independent of $\mdeg$, meaning that the maximum allowable time-step does not degrade as $\mdeg$ increases, which is in contrast to other explicit DG time-stepping approaches (i.e., RKDG and Lax-Wendroff DG).
We observe the expected order of convergence for all methods. 
We also notice that as we consider higher order RIDG methods, the maximum quality increases.
Furthermore, we notice that the $\mdeg=3$ RKDG method has a higher quality than the $\mdeg=3$ RIDG method, however higher order RIDG methods each show a continual increase in the solution quality of 1-2 per 2 increase of $\mdeg$.  Note that each increase of 1 in the solution quality indicates a tenfold decrease in the product of error and runtime (see \cref{eqn:Quality} ).

We can conclude from these results that while the $\mdeg=3$ RKDG method has a superior quality metric than the RIDG method of the same order, the RIDG methods of $\mdeg\geq 5$ have a superior quality metric compared to the RKDG method.  Furthermore, higher order RIDG methods have a superior quality per degree of freedom, since as we increase $\mdeg$ we are able to obtain higher solution qualities with smaller meshes, thus leading to fewer total degrees of freedom.

%% file: tables/advection1d.tex
\begin{table}[!ht]
\centering
\begin{tabular}{|r|r|c|c|c|c|}
\hline	
\multicolumn{6}{|c|}{{\bf 1D RKDG: ($\mdeg=3$, $\nu=0.1$)}}        \\ \hline 
{\bf mesh}&{\bf dof}&{ $L^2$ {\bf error} \eqref{eqn:l2relerror}}&{\bf order} \eqref{eqn:Mratio}&{\bf runtime (s) }&{\bf quality} \eqref{eqn:Quality}\\ \hline %
	$50$	&	$200$	&	$\num{3.07e-4}$	&	--	&	$\num{5.16e-1}$	&	$3.8$	\\
	$70$	&	$280$	&	$\num{7.96e-5}$	&	$4.0$	&	$\num{8.95e-1}$	&	$4.1$	\\
	$120$	&	$480$	&	$\num{9.23e-6}$	&	$4.0$	&	$\num{2.34e+0}$	&	$4.7$	\\
	$240$	&	$960$	&	$\num{5.77e-7}$	&	$4.0$	&	$\num{7.94e+0}$	&	$5.3$	\\ \hline\hline
\multicolumn{6}{|c|}{{\bf 1D RIDG: ($\mdeg=3$, $\nu=0.9$)}}        \\ \hline {\bf mesh}&{\bf dof}&{ $L^2$ {\bf error} \eqref{eqn:l2relerror}}&{\bf order} \eqref{eqn:Mratio}&{\bf runtime (s) }&{\bf quality} \eqref{eqn:Quality}\\ \hline %
	$50$	&	$200$	&	$\num{6.48e-4}$	&	--	&	$\num{1.09e+0}$	&	$3.2$	\\
	$70$	&	$280$	&	$\num{1.47e-4}$	&	$4.4$	&	$\num{1.07e+0}$	&	$3.8$	\\
	$120$	&	$480$	&	$\num{1.85e-5}$	&	$3.8$	&	$\num{2.96e+0}$	&	$4.3$	\\
	$240$	&	$960$	&	$\num{8.31e-7}$	&	$4.5$	&	$\num{1.24e+1}$	&	$5.0$	\\ \hline\hline
\multicolumn{6}{|c|}{{\bf 1D RIDG: ($\mdeg=5$, $\nu=0.9$)}}        \\ \hline {\bf mesh}&{\bf dof}&{ $L^2$ {\bf error} \eqref{eqn:l2relerror}}&{\bf order} \eqref{eqn:Mratio}&{\bf runtime (s) }&{\bf quality} \eqref{eqn:Quality}\\ \hline %
	$30$	&	$180$	&	$\num{1.57e-5}$	&	--	&	$\num{8.29e-1}$	&	$4.9$	\\
	$50$	&	$300$	&	$\num{8.62e-7}$	&	$5.7$	&	$\num{4.21e+0}$	&	$5.4$	\\
	$70$	&	$420$	&	$\num{1.27e-7}$	&	$5.7$	&	$\num{4.17e+0}$	&	$6.3$	\\
	$120$	&	$720$	&	$\num{3.18e-9}$	&	$6.8$	&	$\num{1.24e+1}$	&	$7.4$	\\ \hline\hline
\multicolumn{6}{|c|}{{\bf 1D RIDG: ($\mdeg=7$, $\nu=0.9$)}}        \\ \hline {\bf mesh}&{\bf dof}&{ $L^2$ {\bf error} \eqref{eqn:l2relerror}}&{\bf order} \eqref{eqn:Mratio}&{\bf runtime (s) }&{\bf quality} \eqref{eqn:Quality}\\ \hline %
	$20$	&	$160$	&	$\num{1.32e-6}$	&	--	&	$\num{1.25e+0}$	&	$5.8$	\\
	$30$	&	$240$	&	$\num{4.77e-8}$	&	$8.2$	&	$\num{2.86e+0}$	&	$6.9$	\\
	$50$	&	$400$	&	$\num{6.65e-10}$	&	$8.4$	&	$\num{1.54e+1}$	&	$8.0$	\\
	$70$	&	$560$	&	$\num{6.63e-11}$	&	$6.9$	&	$\num{1.48e+1}$	&	$9.0$	\\ \hline\hline
\multicolumn{6}{|c|}{{\bf 1D RIDG: ($\mdeg=9$, $\nu=0.9$)}}        \\ \hline {\bf mesh}&{\bf dof}&{ $L^2$ {\bf error} \eqref{eqn:l2relerror}}&{\bf order} \eqref{eqn:Mratio}&{\bf runtime (s) }&{\bf quality} \eqref{eqn:Quality}\\ \hline %
	$10$	&	$100$	&	$\num{4.78e-6}$	&	--	&	$\num{1.05e+0}$	&	$5.3$	\\
	$20$	&	$200$	&	$\num{4.62e-9}$	&	$10.0$	&	$\num{4.22e+0}$	&	$7.7$	\\
	$30$	&	$300$	&	$\num{8.25e-11}$	&	$9.9$	&	$\num{9.57e+0}$	&	$9.1$	\\
	$50$	&	$500$	&	$\num{6.11e-13}$	&	$9.6$	&	$\num{4.83e+1}$	&	$10.5$	\\ \hline\hline
\multicolumn{6}{|c|}{{\bf 1D RIDG: ($\mdeg=11$, $\nu=0.9$)}}        \\ \hline {\bf mesh}&{\bf dof}&{ $L^2$ {\bf error} \eqref{eqn:l2relerror}}&{\bf order} \eqref{eqn:Mratio}&{\bf runtime (s) }&{\bf quality} \eqref{eqn:Quality}\\ \hline 
	$5$	&	$60$	&	$\num{1.82e-4}$	&	--	&	$\num{8.72e-1}$	&	$3.8$	\\
	$10$	&	$120$	&	$\num{5.88e-8}$	&	$11.6$	&	$\num{3.52e+0}$	&	$6.7$	\\
	$20$	&	$240$	&	$\num{1.40e-11}$	&	$12.0$	&	$\num{1.28e+1}$	&	$9.7$	\\
	$30$	&	$360$	&	$\num{1.21e-13}$	&	$11.7$	&	$\num{3.08e+1}$	&	$11.4$	\\ \hline
\end{tabular}
\caption{Convergence and runtime study for 1D SSP-RKDG $\mdeg=3$ and RIDG  methods of various $\mdeg$ for the problem defined by \cref{eqn:adv1d}.  We list the mesh size, total number of degrees of freedom, approximate $L^2$ error compared to the exact solution, the approximate observed order of convergence, the wall clock runtime in seconds, and the quality for  each experiment. We see that the RIDG method can use the same CFL $\nu=0.9$ for all values of $\mdeg$. The highest quality solutions are those produced by the highest order RIDG methods. 
\label{table:1D}}
\end{table}%

%% file: section_2d.tex

\subsection{2D results}
\label{sec:2dresults}
Next we consider the RIDG scheme in two spatial dimensions. In particular, we test 
 the full parallel implementation of RIDG by running strong scaling studies to test the manycore capabilities of RIDG as compared to RKDG. 
 Strong scaling is achieved by fixing a problem size and efficiently increasing the amount of compute cores used to solve the same problem. For explicit timestepping methods, this translates to considering a fixed mesh size and subdivide the mesh into smaller pieces, where the number of subdivisions is equal to the number of tasks, as depicted in \cref{table:strong2ddefinition} .  A full strong scaling study subdivides the problem among compute cores until the number of elements (also known as cells or zones) per compute core is minimized (ideally to unity).  This causes the the ratio of MPI communications to work/computation to increase to some maximum value. The 2D strong scaling studies provided in this section are performed at the Institute for Cyber Enabled Research at Michigan State University, using Intel(R) Xeon(R) Gold 6148 2.40GHz CPUs.
\input{tables/strong2ddef.tex}

\subsubsection{2D advection}
For our first test of the scalability of RKDG and RIDG methods in two dimensions, we consider the the 2D scalar advection equation:
\begin{equation}
\label{eqn:adv2d}
\piecewise{    q_{,t} + q_{,x}  + q_{,y} = 0 & (t,\vec{x}) \in \left [ 0, T \right] \times \left [0,1 \right]^2, \\ 
q(t=0,\vec{x}) = e^{\left(\|\vec{x}-\vec{c}\|^2-\omega^2\right)^{-1}} & \text{if} \quad  \|\vec{x}-\vec{c}\|^2 < \omega^2,  \\
q(t=0,\vec{x}) = 0 & \text{if} \quad \|\vec{x}-\vec{c}\|^2 \geq \omega^2, }
\end{equation}
where $\vec{c}=\frac{1}{2}(1,1)^T$ and $\omega = \frac{1}{3}$ and we enforce periodic boundary conditions. We perform the strong scaling study defined by \cref{table:strong2ddefinition} and
report the results in \cref{table:strong2d}.
\input{tables/advection2d_strong.tex}
We notice that the total runtime of the RKDG method is initially less than that of the RIDG method, but does not scale well with granularity. After 1000 cores, the runtime of RKDG increases as we use more compute cores, indicating runtime scaling breakdown.  That is, the many-node RKDG method exhibits low speedups/efficiencies compared to the single node runtime.  The fastest solution was obtained using 360 cores. The RIDG method, while initially more computationally expensive, scales well with the number of added compute resources.  As we add more cores, the RIDG methods of $\mdeg =3,5$ both exhibit monotonic decreases in the overall method runtime, indicating successful runtime scaling to high granularity. As we consider higher granularity, the runtime scaling efficiency for the RIDG method $\mdeg=3$ begins to drop off.  However, the RIDG method $\mdeg=5$ maintains decent runtime scaling efficiency even at very high granularity.  We see that the maximum solution quality obtained by the RKDG method is 5.8, while for RIDG $\mdeg=3$ the maximum quality is 6.5.  For RIDG $\mdeg=5$ the maximum quality is 8.9. Recalling that solution quality is a logarithmic scale, this means that the RIDG $\mdeg=5$ solution exhibits a far superior efficiency.
Although the RIDG method has a very hefty computational cost compared to RKDG, it is able to scale to much higher levels of granularity.  Thus, the RIDG method is able to provide a superior solution quality at scale.

\subsubsection{2D Burgers}
To demonstrate the capability of this method to efficiently solve nonlinear problems, we consider here is the 2D scalar Burgers equation:
\begin{equation}
\label{eqn:burger2d}
\piecewise{    q_{,t} + \paren{\frac{1}{2}q^2}_{,x}  + \paren{\frac{1}{2}q^2}_{,y} = 0 & (t,\vec{x}) \in \left [ 0, T \right] \times \left [0,1 \right]^2, \\ 
q(t=0,\vec{x}) = \frac{1}{4}(1-\cos(x))(1-\cos(y)) & (x,y) \in \left [0,1 \right]^2, }
\end{equation}
where we enforce periodic boundary conditions. We perform the same strong scaling study as described in \cref{table:strong2ddefinition} .  Our results are listed in \cref{table:burger2d} . These studies are also performed at the Institute for Cyber Enabled Research at Michigan State University, using Intel(R) Xeon(R) Gold 6148 2.40GHz CPUs. The  entries of the Jacobian \cref{Jacobian} are calculated using 
the quasi-quadrature-free approach sketched in \cref{sec:jacobians} .
\input{tables/burger2d.tex}

We notice that the RKDG method requires many more total halo-region MPI communications than the RIDG method.  
We also notice that the total runtime of the RKDG method is initially less than that of the RIDG method, but does not scale well with granularity. After 1000 cores, the runtime of RKDG increases as we use more compute cores, indicating runtime scaling breakdown.  The fastest solution was obtained using 360 cores. 
The RIDG method, while initially more computationally expensive, scales well with the number of added compute resources.  As we add more cores, the RIDG method exhibits monotonic decreases in the overall method runtime, indicating successful runtime scaling to high granularity. The fastest runtime was obtained at 1440 cores.
The many-node RKDG method exhibits low speedups/efficiencies compared to the single node runtime.  We notice that at 1000 cores, the speedups drop below $1$ and the efficiencies drop to $0\%$, indicating runtime scaling breakdown. 
The many-node RIDG method exhibits much larger speedups, and thus higher scaling efficiencies. As we consider higher granularity, the runtime scaling efficiency for the RIDG method maintains decent runtime scaling efficiency even at very high granularity.
As we consider various levels of granularity, the maximum solution quality obtained by the RKDG method is 6.0, For RIDG the maximum quality is 7.1. 

We conclude that the RKDG method does not experience runtime scaling to high granularity for this nonlinear problem in two dimensions, but the RIDG method $\mdeg=3$ experiences good runtime scaling to high granularity, providing an improvement in the maximum possible solution quality. This is despite the additional iterations that the RIDG method must perform for nonlinear problems greatly increase the computational cost of the method and the total number of MPI communications.

%% file: tables/strong2ddef.tex
\begin{table}[t!]
\centering
\begin{tabular}{|r|c|c|r|r|}
\hline	
\multicolumn{5}{|c|}{{\bf 2D Strong Scaling Study }}        \\ \hline 
{\bf tasks}&{ \bf mesh/task }&{\bf subdomain}&{\bf total cores}&{\bf elements/core} \\ \hline	
$1$	&	$60^2$	&	$\brak{0.000,1.000}^2$	&	$40$	&	$90.0$	\\
$4$	&	$30^2$	&	$\brak{0.000,0.500}^2$	&	$160$	&	$22.5$	\\
$9$	&	$20^2$	&	$\brak{0.000,0.33\bar{3}}^2$	&	$360$	&	$10.0$	\\
$16$	&	$15^2$	&	$\brak{0.000,0.250}^2$	&	$640$	&	$5.6$	\\
$25$	&	$12^2$	&	$\brak{0.000,0.200}^2$	&	$1000$	&	$3.6$	\\
$36$	&	$10^2$	&	$\brak{0.000,0.16\bar{6}}^2$	&	$1440$	&	$2.5$		\\\hline
\end{tabular}
\caption{Strong scaling study for 2D RKDG and RIDG.  We scale the methods from a fixed physical domain of $\brak{0,1}^2$ to 36 equally sized subdomains.  Each task, which corresponds to one node and thus several processors, solves one subdomain.  Above we present the total number of tasks, the mesh size for each subdivision, the physical size of each subdomain, the total number of cores, and the number of mesh elements per core.  We see that this experiment drives the number of mesh elements per core down to nearly unity. This scaling study was performed at the Institute for Cyber-Enabled Research at Michigan State University, using Intel(R) Xeon(R) Gold 6148 CPU @ 2.40GHz CPUs.  }
\label{table:strong2ddefinition}
\end{table}

%% file: tables/advection2d_strong.tex
\begin{table}[t]
\centering
\begin{tabular}{|r|r|c|c|c|c|c|c|c|c|c|}
\hline	
\multicolumn{8}{|c|}{{\bf 2D RKDG: ($\mdeg=3$, $\nu=0.05$)}}        \\ \hline {\bf cores}&{\bf dof/c}&{\bf efom}&{\bf comms}&{\bf time}&\cref{eqn:Speedup}&\cref{eqn:Efficiency}  &\cref{eqn:Quality}\\ \hline %
	$40$	&	$900.0$	&	$190^2$	&	$9.60{e4}$	&	$3.50{e1}$	&	$--$	&	$--$		&	$5.4$	\\
	$160$	&	$225.0$	&	$95^2$	&	$3.84{e5}$	&	$2.63{e1}$	&	$1.33\times $	&	$11.1\%$	&	$5.6$	\\
	$360$	&	$100.0$	&	$63^2$	&	$8.64{e5}$	&	$1.69{e1}$	&	$2.07\times $	&	$13.4\%$	&	$5.8$	\\
	$1000$	&	$36.0$	&	$38^2$	&	$2.40{e6}$	&	$5.01{e1}$	&	$0.70\times $	&	$0.0\%$		&	$5.3$	\\
	$1440$	&	$25.0$	&	$32^2$	&	$3.46{e6}$	&	$6.99{e1}$	&	$0.50\times $	&	$0.0\%$		&	$5.1$	\\ \hline\hline
\multicolumn{8}{|c|}{{\bf 2D RIDG: ($\mdeg=3$, $\nu=0.7$)}}        \\ \hline {\bf cores}&{\bf dof/c}&{\bf efom}&{\bf comms}&{\bf time}&\cref{eqn:Speedup}&\cref{eqn:Efficiency}  &\cref{eqn:Quality}\\ \hline 
	$40$	&	$1440.0$	&	$240^2$	&	$1.20{e3}$	&	$5.40{e1}$	&	$1.00\times $	&	$--$,	&	$5.5$	\\
	$160$	&	$360.0$	&	$120^2$	&	$4.80{e3}$	&	$1.56{e1}$	&	$3.45\times $	&	$81.7\%$	&	$6.1$	\\
	$360$	&	$160.0$	&	$80^2$	&	$1.08{e4}$	&	$1.37{e1}$	&	$3.94\times $	&	$36.8\%$	&	$6.1$	\\
	$1000$	&	$57.6$	&	$48^2$	&	$3.00{e4}$	&	$6.89{e0}$	&	$7.83\times $	&	$28.5\%$	&	$6.4$	\\
	$1440$	&	$40.0$	&	$40^2$	&	$4.32{e4}$	&	$6.07{e0}$	&	$8.89\times $	&	$22.5\%$	&	$6.5$	\\ \hline\hline
\multicolumn{8}{|c|}{{\bf 2D RIDG: ($\mdeg=5$, $\nu=0.7$)}}        \\ \hline {\bf cores}&{\bf dof/c}&{\bf efom}&{\bf comms}&{\bf time}&\cref{eqn:Speedup}&\cref{eqn:Efficiency}  &\cref{eqn:Quality}\\ \hline 
	$40$	&	$3240.0$	&	$360^2$	&	$1.23{e3}$	&	$8.71{e2}$	&	$1.00\times $	&	$--$	&	$7.6$	\\
	$160$	&	$810.0$	&	$180^2$	&	$4.92{e3}$	&	$2.82{e2}$	&	$3.09\times $	&	$69.5\%$	&	$8.1$	\\
	$360$	&	$360.0$	&	$120^2$	&	$1.11{e4}$	&	$2.37{e2}$	&	$3.68\times $	&	$33.5\%$	&	$8.2$	\\
	$1000$	&	$129.6$	&	$72^2$	&	$3.08{e4}$	&	$7.55{e1}$	&	$11.54\times $	&	$43.9\%$	&	$8.7$	\\
	$1440$	&	$90.0$	&	$60^2$	&	$4.43{e4}$	&	$4.44{e1}$	&	$19.60\times $	&	$53.2\%$	&	$8.9$	\\ \hline
\end{tabular}
\caption{Strong scaling study for the 2D  RIDG  methods with $\mdeg=3$ and $\mdeg=5$ on the 2D advection  \cref{eqn:adv2d}.  We list the total number of cores, the number of degrees of freedom per core, the equivalent first order mesh as described in \cref{eqn:EFOM}, the estimated total number of MPI communications as computed by \cref{eqn:comms}, the runtimes in seconds, the speedup over single node performance as computed by \cref{eqn:Speedup}, the scaling efficiency from a single node \cref{eqn:Efficiency}, and the solution \cref{eqn:Quality}.
\label{table:strong2d}}
\end{table}%

%% file: tables/burger2d.tex
\begin{table}[t]
\centering
\begin{tabular}{|r|r|c|c|c|c|c|c|c|c|c|}
\hline	
\multicolumn{8}{|c|}{{\bf 2D RKDG: ($\mdeg=3$, $\nu=0.05$)}}        \\ \hline {\bf cores}&{\bf dof/c}&{\bf efom}&{\bf comms}&{\bf time}&\cref{eqn:Speedup}&\cref{eqn:Efficiency}  &\cref{eqn:Quality}\\ \hline %
	$40$	&	$900.0$	&	$190^2$	&	$\num{9.60e+4}$	&	$\num{1.90e+1}$	&	$--$	&	$--$	&			$5.5$	\\
	$160$	&	$225.0$	&	$95^2$	&	$\num{3.84e+5}$	&	$\num{7.36e+0}$	&	$2.58\times $	&	$52.7\%$		&	$5.9$	\\
	$360$	&	$100.0$	&	$63^2$	&	$\num{8.64e+5}$	&	$\num{5.51e+0}$	&	$3.45\times $	&	$30.6\%$		&	$6.0$	\\
	$1000$	&	$36.0$	&	$38^2$	&	$\num{2.40e+6}$	&	$\num{1.95e+1}$	&	$0.98\times $	&	$0.0\%$	&			$5.5$	\\
	$1440$	&	$25.0$	&	$32^2$	&	$\num{3.46e+6}$	&	$\num{2.91e+1}$	&	$0.65\times $	&	$0.0\%$	&		$5.3$	\\ \hline\hline
\multicolumn{8}{|c|}{{\bf 2D RIDG: ($\mdeg=3$, $\nu=0.7$)}}        \\ \hline {\bf cores}&{\bf dof/c}&{\bf efom}&{\bf comms}&{\bf time}&\cref{eqn:Speedup}&\cref{eqn:Efficiency}  &\cref{eqn:Quality}\\ \hline 
	$40$	&	$1440.0$	&	$240^2$	&	$\num{5.48e+3}$	&	$\num{8.89e+1}$	&	$1.00\times $	&	$--$		&	$5.9$	\\
	$160$	&	$360.0$	&	$120^2$	&	$\num{2.20e+4}$	&	$\num{2.35e+1}$	&	$3.78\times $	&	$92.6\%$		&	$6.4$	\\
	$360$	&	$160.0$	&	$80^2$	&	$\num{4.92e+4}$	&	$\num{1.26e+1}$	&	$7.07\times $	&	$75.9\%$	&	$6.7$	\\
	$1000$	&	$57.6$	&	$48^2$	&	$\num{1.37e+5}$	&	$\num{7.94e+0}$	&	$11.20\times $	&	$42.5\%$	&	$6.9$	\\
	$1440$	&	$40.0$	&	$40^2$	&	$\num{1.98e+5}$	&	$\num{5.81e+0}$	&	$15.31\times $	&	$40.9\%$		&	$7.1$	\\ \hline
\end{tabular}
\caption{Strong scaling study for Burgers  \cref{eqn:burger2d} using the RIDG methods with $\mdeg=3$ and RKDG methods $\mdeg=3$, as defined by \cref{table:strong2ddefinition}.  We list the total number of cores, the number of degrees of freedom per core, the equivalent first order mesh as described in \cref{eqn:EFOM}, the estimated total number of MPI communications as computed by \cref{eqn:comms}, the runtimes in seconds, the speedup over single node performance as computed by \cref{eqn:Speedup}, the scaling efficiency from a single node \cref{eqn:Efficiency}, and the solution quality \cref{eqn:Quality}.  }
\label{table:burger2d}
\end{table}%

%% file: section_3D_adv.tex

\subsection{3D results}
\label{sec:3dresults}
In this section we test the suitability of the RIDG method at scale for solving 3D problems.
The method efficiency is again compared to the Runge-Kutta DG methods.
In particular we consider the following 3D linear advection equation:
\begin{equation}
\label{eqn:adv3d}
\piecewise{    q_{,t} + q_{,x} + q_{,y} + q_{,z}  = 0 & \text{for} \quad \left(t,\vec{x}\right) \in  \left [ 0, T \right] \times \left [0,1 \right]^3, \\ 
q(t=0,\vec{x}) = \exp{(\|\vec{x}-\vec{c}\|^2-\omega^2 )^{-1}} & \text{if} \quad \|\vec{x}-\vec{c}\|^2 < \omega^2, \\
q(t=0,\vec{x}) = 0  & \text{if} \quad \|\vec{x}-\vec{c}\|^2 \geq \omega^2,} 
\end{equation}
with $\vec{c}=\frac{1}{2}(1,1,1)^T$ and $\omega = \frac{1}{3}$. Note that the initial condition is clearly $C^\infty$ for $\abs{x-c}<\omega$; this is the multidimensional analog of \cref{eqn:adv1d} .

\input{tables/strong3ddef.tex}
\input{tables/advection3d.tex}

To test the scalability of the RIDG method for such a problem, we consider the strong scaling study defined in \cref{table:strong3ddefinition} . As per our scaling studies in the earlier sections, our goal is to shrink the number of degrees of freedom per core.  This maximizes the ratio of communication overhead to task-localized work. In three spatial dimensions, the RIDG prediction stage is a four-dimensional problem and thus there is plenty of task-local work to hide latency, as higher dimensionality exponentially increases the number of elements in each DG basis.  For RIDG, this greatly increases the cost of the region matrix inversions. Despite this our predicted quality metric \cref{eqn:ridgqualityraw} increases as we consider higher order methods. This is demonstrated in the results shown in  \cref{table:advection3D} .  These scaling studies were performed at the Livermore Computing Center at Lawrence Livermore National Laboratory using Intel Xeon E5-2695 v4 2.1GHz CPUs with Omni-Path interconnects. The  entries of the Jacobian \cref{Jacobian} are calculated using 
the quasi-quadrature-free approach sketched in \cref{sec:jacobians}.

First we notice that both the RKDG and RIDG methods exhibit better scaling efficiency for this 3D problem than they do for the 2D problem. The computational complexity of DG increases the ratio of MPI-communication and OpenMP overhead to parallelized work, which also increases the the computational work available to hide the communication latencies, despite the increased number of communications per time-step.  This is indicated by the ``degrees of freedom per core" metric \cref{eqn:ndofspercore} shown in \cref{table:advection3D} compared to the same metric in \cref{table:burger2d} .
Even so, the RKDG method requires many more total halo-region MPI communications than the RIDG method of any order. This is due to the fact that the RIDG methods of $\mdeg=3,5$ exhibit the enhanced CFL restriction $\nu = 0.6$, and thus all can take a relatively large timestep compared to the RKDG method.
 Furthermore, we notice that the total runtime of the RKDG method is half of that of the RIDG method of the same $\mdeg$. Since their errors as per  \cref{table:advection3D} are similar in magnitude, this implies the RKDG method $\mdeg=3$ has a higher quality than the  RIDG method $\mdeg=3$.
 For the higher order RIDG method, $\mdeg=5$, we observe that the runtimes increase by a factor of 2 orders of magnitudes when compared to the methods of $\mdeg=3$.  However, the errors as seen in  \cref{table:strong3ddefinition} are 3 orders of magnitude smaller for $\mdeg=5$. Thus, we find that the higher order methods have a higher quality metric than the methods of $\mdeg=3$, and are thus more efficient.  This is achieved in part by the ability of the RIDG method of $\mdeg=5$ to maintain a strong scaling efficiency of around $90\%$ above 28,000 cores, as the lower order schemes begin to drop off in efficiency.  Note that in 3D the quality metric for ADER-DG methods are approximately $+0.12$ higher than the quality metric for RKDG \cite{Dumbser2018}. 
 We conclude that due to their high quality metric, the high order RIDG method with $\mdeg=5$ scales very efficiently for 3D problems.

%% file: tables/strong3ddef.tex
\begin{table}[t!]
\centering
\begin{tabular}{|r|c|c|r|r|}
\hline	
\multicolumn{5}{|c|}{{\bf 3D Strong Scaling Study }}        \\ \hline 
{\bf tasks}&{ \bf mesh/task }&{\bf subdomain}&{\bf total cores}&{\bf elements/core} \\ \hline	
$1$	&	$48\times 48\times 48$	&	$\brak{0.000, 1.000}^2$	&	$56$	&	$1974.9$	\\
$27$	&	$16\times 16\times 16$	&	$\brak{0.000, 0.33\bar{3}}^2$	&	$1512$	&	$73.1$	\\
$216$	&	$8\times 8\times 8$	&	$\brak{0.000, 0.16\bar{6}}^2$	&	$12096$	&	$9.1$	\\
$512$	&	$6\times 6\times 6$	&	$\brak{0.000, 0.125}^2$	&	$28672$	&	$3.9$	\\
\hline 
\end{tabular}
\caption{Strong scaling study for 3D RKDG and RIDG.  We scale the methods from a fixed physical domain of $\brak{0,1}^3$ to 512 equally sized subdomains.  Each task, which corresponds to one node, solves one subdomain.  Above we present the total number of tasks, the mesh size for each subdivision, the physical size of each subdomain, the total number of cores, and the number of mesh elements per core.  We see that this experiment scales the problem granularity down to about 4 elements per core.  We also display the errors for each method scaled. These runs were performed on Livermore Computing clusters at at Lawrence Livermore National Laboratory using Intel Xeon CLX-8276L CPUs with Omni-Path interconnects. }
\label{table:strong3ddefinition}
\end{table}

%% file: tables/advection3d.tex
\begin{table}[t]
\centering
\begin{tabular}{|r|r|c|c|c|c|c|c|c|c|}
\hline	
\multicolumn{8}{|c|}{{\bf 3D RKDG:} ($\mdeg=3$, $\nu=0.05$, \inred{$L^2 \, \text{\bf error} \, \eqref{eqn:l2relerror} =\num{1.36e-07}$})}        \\ \hline {\bf cores}&{\bf dof/c }&{\bf efom} \cref{eqn:EFOM}&{\bf comms}&{\bf time}&\cref{eqn:Speedup}&\cref{eqn:Efficiency} &\cref{eqn:Quality}\\ \hline %
	$56$	&	$126390.9$	&	$192^3$	&	$\num{2.59e+5}$	&	$\num{4.92e+3}$	&	$--$	&	$--$	&	$4.90$	\\
	$1512$	&	$4681.1$	&	$64^3$	&	$\num{7.00e+6}$	&	$\num{2.15e+2}$	&	$22.86\times $	&	$84.1\%$	&	$6.26$	\\
	$12096$	&	$585.1$	&	$32^3$	&	$\num{5.60e+7}$	&	$\num{2.89e+1}$	&	$170.53\times $	&	$78.8\%$	&	$7.13$	\\
	$28672$	&	$246.9$	&	$24^3$	&	$\num{1.33e+8}$	&	$\num{1.54e+1}$	&	$319.98\times $	&	$62.4\%$	&	$7.41$	\\ \hline\hline
\multicolumn{8}{|c|}{{\bf 3D RIDG:} ($\mdeg=3$, $\nu=0.6$, \inred{$L^2 \, \text{\bf error} \, \eqref{eqn:l2relerror}=\num{1.53e-07}$})}        \\ \hline {\bf cores}&{\bf dof/c }&{\bf efom} \cref{eqn:EFOM}&{\bf comms}&{\bf time}&\cref{eqn:Speedup}&\cref{eqn:Efficiency} &\cref{eqn:Quality}\\ \hline 
	$56$	&	$126390.9$	&	$192^3$	&	$\num{4.32e+3}$	&	$\num{8.24e+3}$	&	$--$	&	$--$	&	$4.72$	\\
	$1512$	&	$4681.1$	&	$64^3$	&	$\num{1.17e+5}$	&	$\num{3.28e+2}$	&	$25.09\times $	&	$92.7\%$	&	$6.13$	\\
	$12096$	&	$585.1$	&	$32^3$	&	$\num{9.33e+5}$	&	$\num{4.43e+1}$	&	$186.17\times $	&	$86.1\%$	&	$7.00$	\\
	$28672$	&	$246.9$	&	$24^3$	&	$\num{2.21e+6}$	&	$\num{2.19e+1}$	&	$377.10\times $	&	$73.6\%$	&	$7.31$	\\ \hline\hline
\multicolumn{8}{|c|}{{\bf 3D RIDG:} ($\mdeg=5$, $\nu=0.6$, \inred{$L^2 \, \text{\bf error} \,\eqref{eqn:l2relerror}=\num{4.97e-11}$})}        \\ \hline {\bf cores}&{\bf dof/c }&{\bf efom} \cref{eqn:EFOM}&{\bf comms}&{\bf time}&\cref{eqn:Speedup}&\cref{eqn:Efficiency} &\cref{eqn:Quality}\\ \hline 
	$1512$	&	$15798.9$	&	$96^3$	&	$1.17{e5}$	&	$2.76{e4}$	&	$--$	&	$--$	&	$7.64$	\\
	$12096$	&	$1974.9$	&	$48^3$	&	$9.33{e5}$	&	$3.76{e3}$	&	$7.34\times $	&	$90.6\%$*	&	$8.51$	\\
	$28672$	&	$833.1$	&	$36^3$	&	$2.21{e6}$	&	$1.57{e3}$	&	$17.54\times $	&	$92.1\%$*	&	$8.88$	\\ \hline
\end{tabular}
\caption{Strong scaling study for the 3D advection \cref{eqn:adv3d} as defined by \cref{table:strong3ddefinition}.  We list the total number of cores, the number of degrees of freedom per core, the equivalent first order mesh as described in \cref{eqn:EFOM}, the estimated total number of MPI communications as computed by \cref{eqn:comms}, the runtimes in seconds, the speedup over single node performance as computed by \cref{eqn:Speedup}, the scaling efficiency from a single node \cref{eqn:Efficiency}, and the solution \cref{eqn:Quality}. In this work we are most interested in the quality metric, as this combines error and runtime into a single metric.  In this table, we see that the high order RIDG methods exhibit the highest quality. Note* the scaling efficiencies \cref{eqn:Efficiency} for the $\mdeg=5$ RIDG scheme are computed relative to the 27 node run as opposed to the single node case. 
} 
\label{table:advection3D}
\end{table}%

%% file: secfinal.tex

\section{Conclusion}
\label{sec:conclusions}
In this paper we have explored the scalability of the regionally implicit discontinuous Galerkin (RIDG) method, which was previously shown to have favorable convergence properties for solving hyperbolic conservation laws \cite{RIDG_paper_2019}.  We compared the results to that of the extremely popular strong-stability-preserving Runge-Kutta DG (SSP-RKDG) method. We performed these comparisons for a two-dimensional linear problem, a two-dimensional nonlinear problem, and a three-dimensional linear problem which served as a toy problem related to the Relativistic Vlasov Maxwell system. We demonstrated efficient strong-scaling properties of the RIDG DG methods, we are able to maintain vertex-neighbor stenciling and communications while taking highly enhanced time-step sizes.  We also demonstrated the incredible boost in efficiency of using the quasi-quadrature free strategy for space-time matrix assembly.  This strategy reduces the computational complexity of space-time matrix assembly down to what we would expect from purely spatial matrix-assembly. This is critical to the viability of RIDG methods.

    The work used to hide the aforementioned communication and collective latencies itself benefits from intra-node (shared memory) parellelism via OpenMP \cite{openmp,Architecture2018}. The RIDG method achieves decent intra-node efficiency simply by distributing the ``regions'' over the number of shared memory compute cores, similar to the ADER-DG parellelism strategy \cite{Fambri2018}. The result is that in the limit of minimum granularity (minimum cells per core), the work done by each core is dominated by formation of the space-time region Jacobian and small linear system solves.  In contrast, the RKDG method distributes the computation of the flux quadrature and the volume integral quadrature to form a residual.  In the the limit of minimum granularity, this means that each CPU core computes the spatial component of the residual of a cell.  Thus, the RIDG method distributes a far greater load of work compared to the RKDG method, which leads to superior intra-node scaling, at the cost of incredibly increased computational cost per time step.
   
   Further, the superior stability properties of the RIDG method allowed larger time steps than what the RKDG method was able to take. The predictor-corrector strategy of the RIDG method means that the method communicates twice per time-step compared to once per stage (many times per time-step).  Thus, the RIDG method communicates across domain decomposition pseudo-boundaries (halo regions) nearly two orders of magnitudes fewer times than the RKDG method.  In addition, the amount of work used to hide the latency of the communications is much greater for the RIDG method than for the RKDG method.  In the nonlinear case, this advantage is diminished as the RIDG method must communicate once per Newton iteration during the prediction step, but the RIDG method still exhibits a superior ability to hide communication latency behind the Newton iteration matrix solves. 
    
    Furthermore, the reduced number of time-steps means that there are fewer all-to-all collective communications of the maximum wavespeed, an operation required of explicit time integrators. The RIDG method is able to compute the maximum wavespeed after the prediction step, and is thus able to hide the latency associated with this collective using the corrector step. 
    
    Lastly, reduction in parallel scaling efficiency is inevitable, as a growing percentage of each subdomain must be copied into contiguous MPI buffers and a growing number of compute resources must be synchronized for maximum efficiency.  However, if the ratio of time spent updating the solution over the time spent performing these copies is maximized, as is with RIDG, then higher parallel scaling can be maintainted.
    
Due to excellent inter-node and intra-node parallelism, the RIDG method offers excellent scalability for explicit time-stepping of hyperbolic conservation laws, and is able to demonstrate incredible performance with very high-order methods at high dimensionality.

%% file: main_bit.bbl
\begin{thebibliography}{10}
\providecommand{\url}[1]{{#1}}
\providecommand{\urlprefix}{URL }
\expandafter\ifx\csname urlstyle\endcsname\relax
  \providecommand{\doi}[1]{DOI~\discretionary{}{}{}#1}\else
  \providecommand{\doi}{DOI~\discretionary{}{}{}\begingroup
  \urlstyle{rm}\Url}\fi

\bibitem{Atkins1998}
Atkins, H.L., Shu, C.W.: Quadrature-free implementation of discontinuous
  {G}alerkin method for hyperbolic equations.
\newblock AIAA Journal \textbf{36}(5), 775--782 (1998).
\newblock \doi{10.2514/3.13891}.
\newblock \urlprefix\url{http://arc.aiaa.org/doi/10.2514/2.436}

\bibitem{chan2016}
Chan, J., Wang, Z., Modave, A., Remacle, J.F., Warburton, T.: {GPU-accelerated
  discontinuous Galerkin methods on hybrid meshes}.
\newblock J. Comput. Physics \textbf{318}(div), 142--168 (2016).
\newblock \doi{10.1016/j.jcp.2016.04.003}

\bibitem{cockshu5}
Cockburn, B., Shu, C.W.: The {R}unge--{K}utta discontinuous {G}alerkin method
  for conservation laws {V}: {M}ultidimensional systems.
\newblock J. Comput. Physics \textbf{141}(2), 199--224 (1998)

\bibitem{openmp}
{Dagum, Leonardo and Menon}, R.: {OpenMP: an industry standard API for
  shared-memory programming}.
\newblock Computational Science \& Engineering, IEEE \textbf{5}, 46--55 (1998)

\bibitem{Dumbser2018}
Dumbser, M., Fambri, F., Tavelli, M., Bader, M., Weinzierl, T.: {Efficient
  implementation of ADER discontinuous Galerkin schemes for a scalable
  hyperbolic PDE engine}.
\newblock Axioms pp. 1--26 (2018).
\newblock \doi{10.3390/axioms7030063}.
\newblock \urlprefix\url{http://arxiv.org/abs/1808.03788}

\bibitem{article:Dumbser2006}
Dumbser, M., Munz, C.D.: Building blocks for arbitrary high order discontinuous
  {G}alerkin schemes.
\newblock J. Sci. Comput. \textbf{27}, 215--230 (2006)

\bibitem{Dumbser2013a}
Dumbser, M., Zanotti, O., Hidalgo, A., Balsara, D.S.: {ADER-WENO finite volume
  schemes with space-time adaptive mesh refinement}.
\newblock J. Comput. Physics \textbf{248}, 257--286 (2013).
\newblock \doi{10.1016/j.jcp.2013.04.017}.
\newblock \urlprefix\url{http://dx.doi.org/10.1016/j.jcp.2013.04.017}

\bibitem{Engsig-Karup2016}
Engsig-Karup, A.P., Eskilsson, C., Bigoni, D.: {A stabilised nodal spectral
  element method for fully nonlinear water waves}.
\newblock J. Comput. Physics \textbf{318}, 1--21 (2016).
\newblock \doi{10.1016/j.jcp.2016.04.060}

\bibitem{Fambri2018}
Fambri, F., Dumbser, M., K{\"{o}}ppel, S., Rezzolla, L., Zanotti, O.: {ADER
  discontinuous Galerkin schemes for general-relativistic ideal
  magnetohydrodynamics}.
\newblock Mon. Not. R. Astron. Soc. \textbf{477}(4), 4543--4564 (2018).
\newblock \doi{10.1093/mnras/sty734}

\bibitem{Fischer}
Fischer, P.F., Heisey, K., Min, M.: {Scaling Limits for PDE-Based Simulation}.
\newblock In: 22nd AIAA Computational Fluid Dynamics Conference, 22-26 June
  2015, Dallas, TX, pp. 1--10 (2015)

\bibitem{Gabriel2004}
Gabriel, E., Fagg, G.E., Bosilca, G., Angskun, T., Dongarra, J.J., Squyres,
  J.M., Sahay, V., Kambadur, P., Barrett, B., Lumsdaine, A., Castain, R.H.,
  Daniel, D.J., Graham, R.L., Woodall, T.S.: {Open MPI: Goals, concept, and
  design of a next generation MPI implementation}.
\newblock Lecture Notes in Computer Science (including subseries Lecture Notes
  in Artificial Intelligence and Lecture Notes in Bioinformatics)
  \textbf{3241}, 97--104 (2004).
\newblock \doi{10.1007/978-3-540-30218-6-19}

\bibitem{article:GasDumHinMun2011}
Gassner, G., Dumbser, M., Hindenlang, F., Munz, C.D.: Explicit one-step time
  discretizations for discontinuous {G}alerkin and finite volume schemes based
  on local predictors.
\newblock J. Comput. Physics \textbf{230}, 4232--4247 (2011)

\bibitem{article:GoShu98}
Gottlieb, S., Shu, C.W.: Total variation diminshing {R}unge-{K}utta schemes.
\newblock Math. of Comput. \textbf{67}, 73--85 (1998)

\bibitem{gottliebShuTadmor01}
Gottlieb, S., Shu, C.W., Tadmor, E.: Strong stability-preserving high-order
  time discretization methods.
\newblock SIAM Rev. \textbf{43}(1), 89--112 (2001)

\bibitem{Gupta2007}
Gupta, M., Narasimhan, S.G.: {Legendre polynomials Triple Product Integral and
  lower-degree approximation of polynomials using Chebyshev polynomials}.
\newblock Tech. Rep. CMU-RI-TR-07-22, Carnegie Mellon University (2007)

\bibitem{Guthrey2017}
Guthrey, P.: {Regionally implicit discontinuous Galerkin methods for solving
  the relativistic Vlasov-Maxwell system}.
\newblock Ph.D. thesis, Iowa State University (2017)

\bibitem{RIDG_paper_2019}
Guthrey, P., Rossmanith, J.: {The regionally implicit discontinuous Galerkin
  method: Improving the stability of DG-FEM}.
\newblock SIAM J. Numer. Analysis \textbf{57}(3), 1263--1288 (2019).
\newblock \doi{10.1137/17M1156174}

\bibitem{Consortium2017}
{INTERTWinE Consortium}: {Best Practice Guide to Hybrid MPI + OpenMP
  Programming} (2017).
\newblock
  \urlprefix\url{http://www.intertwine-project.eu/sites/default/files/images/INTERTWinE_Best_Practice_Guide_MPI%2BOpenMP_1.1.pdf}

\bibitem{article:KlaVegVen2006}
Klaij, C., van Der~Vegt, J., van Der~Ven, H.: Space-time discontinuous
  {G}alerkin method for the compressible {N}avier-{S}tokes equations.
\newblock J. Comput. Physics \textbf{217}(2), 589--611 (2006)

\bibitem{Architecture2018}
OpenMP Architecture Review Board: {OpenMP Application Programming Interface
  } (2018).
\newblock
  \urlprefix\url{https://www.openmp.org/wp-content/uploads/OpenMP-API-Specification-5.0.pdf}

\bibitem{Qiu2005a}
Qiu, J., Dumbser, M., Shu, C.W.: {The discontinuous Galerkin method with
  Lax-Wendroff type time discretizations}.
\newblock Comput. Methods Appl. Mech. Eng. \textbf{194}, 4528--4543 (2005).
\newblock \doi{10.1016/j.cma.2004.11.007}

\bibitem{article:Ru61}
Rusanov, V.: Calculation of interaction of non-steady shock waves with
  obstacles.
\newblock J. Comp. Math. Phys. USSR \textbf{1}, 267--279 (1961)

\bibitem{article:Sudirham2006}
Sudirham, J., van Der~Vegt, J., van Damme, R.: Space-time discontinuous
  {G}alerkin method for advection-diffusion problems on time-dependent domains.
\newblock Appl. Numer. Math. \textbf{56}(12), 1491--1518 (2006)

\bibitem{article:Zanotti2015}
Zanotti, O., Fambri, F., Dumbser, M., Hidalgo, A.: {Space--time adaptive ADER
  discontinuous Galerkin finite element schemes with a posteriori sub-cell
  finite volume limiting}.
\newblock Computers {\&} Fluids \textbf{118}, 204 -- 224 (2015).
\newblock \doi{https://doi.org/10.1016/j.compfluid.2015.06.020}

\end{thebibliography}
